\documentclass{amsart}
\usepackage{amssymb,amsthm}
\usepackage[arrow,matrix]{xy}
\usepackage{pstricks}

\newcommand{\nc}{\newcommand}
\nc{\rnc}{\renewcommand}

\nc{\al}{\alpha}
\nc{\be}{\beta}
\nc{\ga}{\gamma}
\nc{\de}{\delta}
\nc{\eps}{\varepsilon}
\nc{\la}{\lambda}
\nc{\sig}{\sigma}
\nc{\om}{\varpi}

\nc{\<}{\langle}
\rnc{\>}{\rangle}

\nc{\xto}{\xrightarrow}

\nc{\C}{\mathbf C}
\nc{\N}{\mathbf N}

\rnc{\<}{\langle}
\rnc{\>}{\rangle}

\rnc{\P}{\mathbb P}
\nc{\MP}{\P^{1\dots\ell}}
\nc{\MT}{\Tens(V_1\oplus\dots\oplus V_\ell)}
\nc{\g}{\mathfrak g}
\rnc{\sl}{\mathfrak{sl}}
\nc{\hi}{\widehat\imath}
\nc{\DrJ}{\mathrm{DJ}}

\nc{\pr}{\operatorname{pr}}
\nc{\Pl}{\operatorname{Pl}}
\nc{\Norm}{\operatorname{N}}
\nc{\Triple}{\operatorname{T}}
\nc{\Tens}{\operatorname{T}}

\nc{\D}{\Delta}
\rnc{\L}{{\mathcal L}}
\rnc{\O}{\operatorname{\mathcal O}}
\nc{\HZ}{\operatorname H^0}

\nc{\simto}{\overset{\sim}\to}
\nc{\from}{\leftarrow}
\rnc{\o}{\otimes}
\nc{\ol}{\overline}
\nc{\inv}{^{-1}}
\nc{\Pp}{P^{+}}
\nc{\rien}{\text{--}}
\nc{\exch}{\leftrightarrow}
\nc{\spmat}[1]{\left(\begin{smallmatrix}#1\end{smallmatrix}\right)}
\rnc{\labelenumi}{(\alph{enumi})}
\nc{\sll}{\operatorname{\mathfrak{sl}}}
\nc{\SL}{\operatorname{SL}}
\nc{\GL}{\operatorname{GL}}
\nc{\SO}{\operatorname{SO}}
\rnc{\Sp}{\operatorname{Sp}}
\nc{\U}{\operatorname{U}}
\nc{\Uqsln}{\U_q^\DrJ(\sll(n))}
\nc{\Uqg}{\U_q^\DrJ(\g)}
\nc{\OqG}{\mathcal O_q^\DrJ(G)}
\nc{\OqSLn}{\mathcal O_q^\DrJ(\SL(n))}
\nc{\Segre}{\text{Segre}}
\nc{\diag}{\text{diag.}}

\nc{\id}{\mathrm{id}}

\theoremstyle{definition}
\newtheorem{df}{Definition}[section]
\newtheorem{rem}[df]{Remark}
\newtheorem{quest}{Question}
\newtheorem{prob}[quest]{Problem}
\newtheorem{ex}[df]{Example}

\theoremstyle{plain}
\newtheorem{theo}[df]{Theorem}
\newtheorem{conj}[df]{Conjecture}
\newtheorem{prop}[df]{Proposition}
\newtheorem{lem}[df]{Lemma}
\newtheorem{cor}[df]{Corollary}
\newtheorem{crit}[df]{Criterion}

\numberwithin{equation}{section}

\rnc{\thequest}{\Alph{quest}}
\rnc{\theprob}{\Alph{prob}}

\rnc{\subjclassname}{%
  \textup{2000} Mathematics Subject Classification}

\title[A ``classical'' flag variety for standard quantum
$\SL(n)$]{``Classical'' flag varieties for quantum groups: the standard
quantum $\SL(n,\C)$}
\author{Christian Ohn}
\date{June 30, 2000; this version January 12, 2001}
\address{Universit\'e de Reims
\\
D\'epartement de Math\'ematiques (UPRESA 6056 du CNRS)
\\
Moulin de la Housse, B.P. 1039
\\
F-51687 Reims Cedex 2
\\
France}
\email{christian.ohn@univ-reims.fr}

\begin{document}

\begin{abstract}
We suggest a possible programme to associate geometric ``flag-like''
data to an arbitrary simple quantum group, in the spirit of the
noncommutative algebraic geometry developed by Artin, Tate, and
Van~den~Bergh. We then carry out this programme for the standard
quantum $\SL(n)$ of Drinfel$'$d and Jimbo, where the varieties involved
are certain $T$-stable subvarieties of the (ordinary) flag variety.
\end{abstract}

\maketitle

\part*{Introduction}

The study of quantum analogues of flag varieties, first suggested by
Manin~\cite{Ma}, has been undertaken during the past decade by several
authors, from various points of view; see e.g.\
\cite{APW,Br,DG,Fi,Go,Jo,LaRe,LuRo,So,SD,ST,TT}. Around the same time, an
approach to noncommutative projective algebraic geometry was initiated by
Artin, Tate, and Van~den~Bergh~\cite{ATV1,ATV2,AV}, and considerably
developed since (see e.g.~\cite{AS,AZ,Ca,K,LSV,Stp,VV,VdB,Z}). One attractive
feature of their approach is the association of actual \emph{geometric} data
to certain classes of graded \emph{non}commutative algebras.

The present work is an attempt to study quantum flag varieties from this point
of view. As a consequence, our ``quantum flag varieties'' will be actual
varieties (with some bells and whistles).

Recall the original idea of~\cite{ATV1,ATV2,AV}. If $A$ is the homogeneous
coordinate ring of a projective scheme $E$, then the points of $E$ are in
one-to-one correspondence with the isomorphism classes of so-called
\emph{point modules} of $A$, i.e.\ $\N$-graded cyclic $A$-mod\-ules $P$ such
that $\dim P_n=1$ for all $n$. Now if $A$ is an $\N$-graded
\emph{non}commutative algebra, one may still try to parametrize its point
modules by the points of some projective scheme $E$. Of course, one cannot
hope to reconstruct $A$ from $E$ alone, but there is now an additional
ingredient: the \emph{shift} operation $\sig:P\mapsto P[1]$, where $P[1]$ is
the $\N$-graded $A$-mod\-ule defined by $P[1]_n:=P_{n+1}$. (When $A$ is
commutative, this shift is trivial: $P[1]\simeq P$ for every point module
$P$.) Assume that $\sig$ may be viewed as an automorphism of $E$: one may
then hope, at least in ``good'' cases, to recover $A$ from the triple
$(E,\sig,\L)$, where $\L$ is the line bundle over $E$ defined by its
embedding into a projective space. The first step of this recovery is the
construction of the \emph{twisted homogeneous coordinate ring} $B(E,\sig,\L)$
of a triple $(E,\sig,\L)$, defined in~\cite{AV} as follows:
\[
B(E,\sig,\L)=\bigoplus_{n\in\N}B_n,\qquad
B_n:=\HZ(E,\L\o\L^\sig\o\dots\o\L^{\sig^{n-1}})
\]
(where $\L^\sig$ denotes the pullback of $\L$ along $\sig$), the
multiplication being given by $\al\be:=\al\o\be^{\sig^m}$ for all $\al\in
B_m$, $\be\in B_n$. (When $\sig$ is the identity, this algebra is the
(commutative) homogeneous coordinate ring of $E$ w.r.t.\ the polarization
$\L$.) If the triple $(E,\sig,\L)$ comes from an algebra $A$ as above, the
second step then consists in analysing the canonical morphism $A\to
B(E,\sig,\L)$. The initial success of this method has been a complete study
of all regular algebras of dimension three~\cite{ATV1} (where the kernel of
$A\to B(E,\sig,\L)$ turns out to be generated by a single element of degree
three).

The present paper is organized as follows.

In Part~\ref{outlinepart}, we give a general outline of a possible theory of
flag varieties for quantum groups, using a multigraded version of the ideas
of~\cite{ATV1,ATV2,AV} recalled above, some of which have already been
introduced by Chan~\cite{Ch}. This Part is largely conjectural and contains
no (significant) new results; its purpose is rather to set up a framework that
will be used in Parts \ref{djpart}~and \ref{djslnpart}.

More specifically, we proceed as follows. Let $G$ be a simple complex group;
our interest in flag varieties allows us to assume without harm that $G$ is
simply connected. Let $\Pp$ be the monoid of dominant integral weights of $G$
(w.r.t.\ some Borel subgroup $B\subset G$): the \emph{shape algebra} $M$ of
$G$ is a $\Pp$-graded $G$-alg\-ebra whose term of degree $\la$ is the
irreducible representation of $G$ of highest weight $\la$. Now consider the
definition of a point module (see above), but with $\N$-grad\-ings replaced
by $\Pp$-grad\-ings: we obtain the notion of a \emph{flag module} of $M$
(Definition~\ref{flagmoddef}); this terminology is justified by the fact that
the isomorphism classes of such modules are indeed parametrized by the points
of the flag variety $G/B$ (Proposition~\ref{commflags}).

If a quantum group has the ``same'' representation theory as $G$ (in the
sense of Definition~\ref{Gqdef}), then the we may still define a
($\Pp$-graded) shape algebra. We then discuss the possibility to parametrize
the latter's flag modules (up to isomorphism) by the points of some scheme
$E$, and to realize shifts $F\mapsto F[\la]$ ($\la\in\Pp$) as automorphisms
$\sig_\la$ of $E$. It will of course be sufficient to know the automorphisms
$\sig_1,\dots,\sig_\ell$ associated to the fundamental weights
$\om_1,\dots,\om_\ell$, which freely generate $\Pp$. Moreover, since we are
in a \emph{multi}graded situation, it will be more natural to view $E$ as a
subscheme of a \emph{product} of $\ell$ projective spaces, corresponding to
$\ell$ line bundles $\L_1,\dots,\L_\ell$ over $E$.

We then consider the converse problem of reconstructing the shape algebra
from $E$, the $\sig_i$, and the $\L_i$, using Chan's construction~\cite{Ch}
of a twisted \emph{multi}homogeneous coordinate ring: this is the
$\Pp$-graded algebra
\[
B(E,\sig_1,\dots,\sig_\ell,\L_1,\dots,\L_\ell)
:=\bigoplus_{\la\in\Pp}\HZ(E,\L_\la),
\]
where the line bundles $\L_\la$ are constructed inductively from the rules
$\L_{\om_i}=\L_i$, $\L_{\om_i+\la}=\L_i\o\L_\la^{\sig_i}$. (Again, if
$E=G/B$, if each $\sig_i$ is the identity, and if the $\L_i$ are the line
bundles associated to the fundamental $G$-mod\-ules $V^1,\dots,V^\ell$, then
this algebra is the (commutative) multihomogeneous coordinate ring of
$G/B\subset\P(V^1)\times\dots\times\P(V^\ell)$, which in turn is equal to the
shape algebra $\O(\ol{G/U})$, $U$ the unipotent radical of $B$.)

We stress that the ideas developed in this Part are \emph{not} restricted to
the standard quantum groups of Drinfel$'$d~\cite{Dr} and Jimbo~\cite{Ji}, but
could, in principle, be applied to other quantum groups as well, as long as
they have the ``same'' representation theory as a given simple complex group.
(Potential other examples include the multiparameter quantum groups of Artin,
Schelter, and Tate~\cite{AST,HLT}, the quantum $\SL(n)$ of Cremmer and
Gervais~\cite{CG,Ho}, or the quantum $\SL(3)$'s classified in~\cite{qsl3}.)

In Parts \ref{djpart}~and \ref{djslnpart}, we \emph{do} restrict ourselves to a
standard Drinfel$'$d-Jimbo quantum group $\OqG$, with $q$ not a root of
unity. Thanks to the results of Lusztig~\cite{Lu} and Rosso~\cite{Ro}, $\OqG$
has the ``same'' representation theory as the group $G$, so one can define a
shape algebra $M^\DrJ$.

In Part~\ref{djpart}, we construct geometric data $E^\DrJ$, $\sig_i$, and
$\L_i$, and we conjecture that these data indeed correspond to $M^\DrJ$ as
described above (Conjectures \ref{mainconj}~and \ref{secconj}). The scheme
$E^\DrJ$ will actually be a union of certain $T$-stable subvarieties of the
(ordinary) flag variety $G/B$. Since the latter may be of independent
interest to algebraic geometers, we have decided to describe them in a
separate note~\cite{PPinGB} (but we recall their construction here, without
proofs).

In Part~\ref{djslnpart}, we prove Conjecture~\ref{mainconj} for $G=\SL(n)$,
thus obtaining a ``flag variety'' for the standard Drinfel$'$d-Jimbo quantum
$\SL(n)$. The proof uses special features of the group $\SL(n)$ (the Weyl
group is the symmetric group, all fundamental representations are minuscule,
etc.) and is essentially combinatorial; it is therefore not likely to be
extendable to an arbitrary $G$.

{\bf Acknowledgement.} The author would like to thank the Universit\'e
de Reims for granting a sabbatical leave during the year 1999--2000,
when part of this work has been done.

{\bf Conventions.} All vector spaces, dimensions, algebras, tensor products,
varieties, schemes, etc.\ will be over the field $\C$ of complex numbers. If
$G$ is a linear algebraic group, we denote by $\O(G)$ the Hopf algebra of
polynomial functions on $G$. If $A$ is a (co)algebra, then the dual of a left
$A$-(co)mod\-ule is a right $A$-(co)module, and vice-versa; morphisms of
$A$-(co)mod\-ules will simply be called \emph{$A$-morph\-isms}. When $V$ is a
vector space and $v\in V$ is nonzero, we will sometimes still denote by $v$
the corresponding point in the projective space $\P(V)$.

\tableofcontents

\part{An approach to quantum flag varieties: general outline}
\label{outlinepart}

This Part contains no (significant) new results. It rather discusses a
possible theory of flag varieties for simple quantum groups, asking several
questions along the way (as well as two ambitious problems at the end).

Most of what we will say here is a multigraded version of some of the
main ideas of~\cite{ATV1,ATV2,AV}, applied in a Lie-theoretic setting.

\section{Simple quantum groups and their shape algebras}
\label{shapesection}

Let $G$ be a simply connected simple complex group, $B\subset G$ a Borel
subgroup, and $\Pp$ the set of dominant integral weights of $G$
w.r.t.\ $B$. For each $\la,\mu,\nu\in\Pp$, denote by
\begin{itemize}
\item $d^\la$ the dimension of the simple $G$-mod\-ule of highest weight
$\la$, and by
\item $c^{\la\mu}_\nu$ the multiplicity of the simple $G$-mod\-ule of
highest weight $\nu$ inside the tensor product of those of highest weights
$\la$ and $\mu$.
\end{itemize}
Bearing in mind that the algebra $\O(G)$ of polynomial functions on $G$
is a commutative Hopf algebra, and that (finite-dimensional) left
$G$-mod\-ules correspond to right $\O(G)$-co\-mod\-ules, recall the following
definition from~\cite{qsl3}.
\begin{df}
\label{Gqdef}
We call a \emph{quantum $G$} any (not necessarily commutative) Hopf
algebra $A$ (over $\C$) such that
\begin{enumerate}
\item there is a family $\{V^\la\mid\la\in\Pp\}$ of simple and pairwise
nonisomorphic (right) $A$-co\-mod\-ules, with $\dim V^\la=d^\la$,
\item every $A$-co\-mod\-ule is isomorphic to a direct sum of these,
\item for every $\la,\mu\in\Pp$, $V^\la\o V^\mu$ is isomorphic to
$\bigoplus_{\nu} c^{\la\mu}_\nu V^\nu$.
\end{enumerate}
\end{df}
For convenience, we will write
\[
V_\la:=(V^\la)^*.
\]
For every $\la,\mu\in\Pp$, Definition~\ref{Gqdef}(c) yields an injective
$A$-morph\-ism $V^{\la+\mu}\to V^\la\o V^\mu$ that is unique up to scalars.
Denote by
\[
m_{\la\mu}:V_\la\o V_\mu\to V_{\la+\mu}
\]
the corresponding projection. Gluing these together on
\[
M_A:=\bigoplus_\la V_\la,
\]
we get a (not necessarily associative) multiplication $m:M_A\o M_A\to M_A$.
\begin{df}
The algebra $M_A$ is called the \emph{shape algebra} of $A$.
\end{df}
\begin{quest}
\label{assoc}
Is it possible to renormalize the $m_{\la\mu}$ in such a way that the
multiplication $m$ becomes associative?
\end{quest}
Recall that this Question has a positive answer in the commutative case
$A=\O(G)$: if $U$ is a maximal unipotent subgroup, then by the Borel-Weil
theorem, we may set
\[
M_{\O(G)}=\O(\ol{G/U}):=\{f\in\O(G)\mid f(gu)=f(g)\ \forall g\in G,\ \forall
u\in U\}.
\]
The next Proposition provides a criterion for a positive answer to
Question~\ref{assoc}. We first introduce some more notation:
let $\ell$ be the rank of $G$, denote by $\om_1,\dots,\om_\ell$ the
fundamental weights, and let us use the shorthand notation
\[
V_i:=V_{\om_i},\qquad V^i:=V^{\om_i},\qquad1\le i\le\ell.
\]
For every $1\le i,j,k\le\ell$, Definition~\ref{Gqdef}(c) implies that $V_i\o
V_j\o V_k$ contains a unique subcomodule isomorphic to
$V_{\om_i+\om_j+\om_k}$; denote this subcomodule by $W_{ijk}$.
\begin{prop}
\label{assocbraid}
Question~\emph{\ref{assoc}} has a positive answer (for a given $A$) if and
only if there exist $A$-iso\-morph\-isms $R_{ij}:V_i\o V_j\to V_j\o V_i$ for
all $i>j$, such that the braid relation
\begin{equation}
\label{Vijkbraid}
(R_{jk}\o\id)(\id\o R_{ik})(R_{ij}\o\id)\mid_{W_{ijk}}
=(\id\o R_{ij})(R_{ik}\o\id)(\id\o R_{jk})\mid_{W_{ijk}}
\end{equation}
holds for all $i>j>k$.
\end{prop}
\noindent We defer the proof to Appendix~\ref{assocbraidproof}.
\begin{cor}
Question~\emph{\ref{assoc}} has a positive answer in
each of the following situations:
\begin{itemize}
\item when $G$ is of rank~$2$,
\item when $A$ is dual quasitriangular,
\item when $G=\SL(n)$ (by the main result of~\cite{KW}).
\end{itemize}
\end{cor}
Since $\om_1,\dots,\om_\ell$ generate the monoid $\Pp$, and since the
$m_{\la\mu}$ are surjective, the algebra $M_A$ is generated by
\[
M_1:=V_1\oplus\dots\oplus V_\ell.
\]
In this way, $M_A$ may be viewed as an $\N$-graded algebra. More explicitly,
if $\la\in\Pp$ decomposes as
$\sum_i a_i\om_i$ (with each $a_i\in\N$), and if we write $h(\la):=\sum a_i$
for the \emph{height} of $\la$, then the $\N$-grad\-ing on $M_A$ is given by
$M_k:=\bigoplus_{h(\la)=k}V_\la$.
\begin{quest}
\label{quadratic}
Is the shape algebra $M_A$ quadratic (as an $\N$-graded algebra)?
\end{quest}
In the commutative case $A=\O(G)$, the shape algebra $\O(\ol{G/U})$ is indeed
quadratic by a well known theorem of Kostant (see~\cite[Theorem~1.1]{LT} for
a proof). This remains true for the standard Drinfel$'$d-Jimbo quantum
$\SL(n)$: a presentation of the corresponding shape algebra by generators and
(quadratic) relations has been given by Taft and Towber~\cite{TT}.
\begin{quest}
Is the shape algebra $M_A$ a Koszul algebra?
\end{quest}
To finish this Section, let us take a closer look at the quadratic relations
in $M_A$. For every $1\le i,j\le\ell$, let $K_{ij}$ be the kernel of the
multiplication $V_i\o V_j\to V_{\om_i+\om_j}$. By Definition~\ref{Gqdef}(c),
the $A$-co\-mod\-ules $V_i\o V_j$ and $V_j\o V_i$ are isomorphic, and,
rescaling the $A$-iso\-morph\-ism $R_{ij}:V_i\o V_j\to V_j\o V_i$ of
Proposition~\ref{assocbraid} if necessary, we may assume that the diagram
\eqref{scaleRij} (in Appendix~\ref{assocbraidproof}) commutes. Using
Definition~\ref{Gqdef}(c), we see that the quadratic relations in $M_A$ of
degree $\om_i+\om_j$ are of two kinds:
\begin{enumerate}
\rnc{\labelenumi}{(\Roman{enumi})$_{ij}$}
\item $\xi=0$, for $\xi\in K_{ij}$;
\item $\xi=R_{ij}(\xi)$, for $\xi\in V_i\o V_j$.
\end{enumerate}
\begin{rem}
\label{superfluous}
By Definition~\ref{Gqdef}(c), relations (I)$_{ij}$ and (II)$_{ij}$ for
arbitrary $i,j$ are consequences of relations (I)$_{ij}$ for $i\ge j$ only and
relations (II)$_{ij}$ for $i>j$ only.
\end{rem}

\section{The scheme of flag modules}
\label{schemeflagsection}

Assume that Question~\ref{assoc} has a positive answer. The following
definition is a multigraded analogue of the point modules introduced
in~\cite{ATV2}.
\begin{df}
\label{flagmoddef}
A \emph{flag module} is a $\Pp$-graded right
$M_A$-mod\-ule $F$ such that
\begin{enumerate}
\item $\dim F_\la=1$ for each $\la\in\Pp$,
\item $F$ is cyclic.
\end{enumerate}
\end{df}
The terminology is justified by the commutative case. Indeed, let $B\subset G$
be a Borel subgroup and $U$ the unipotent radical of $B$. Then we have the
following
\begin{prop}
\label{commflags}
The isomorphism classes of flag modules of $M_{\O(G)}=\O(\ol{G/U})$ are
pa\-ra\-metrized by the points of the flag variety $G/B$.
\end{prop}
\begin{proof}
First, recall from the Borel-Weil theorem that the decomposition
$\O(\ol{G/U})=\bigoplus_{\la\in\Pp}V_\la$ is given by
\[
V_\la=\{f\in\O(G)\mid f(gb)=\la(b)f(g)\ \forall g\in G,\ \forall b\in B\}.
\]
Now fix $g\in G$ and endow a vector space $F=\bigoplus_{\la\in\Pp}\C e_\la$
with the flag module structure defined by
\[
e_\la.f=f(g)e_{\la+\mu}\qquad\text{for all $f\in V_\mu$.}
\]
If we replace $g$ by $gb$ for some $b\in B$, the expression for
$e_\la.f$ is just multiplied by $\mu(b)$, so up to isomorphism (of
graded modules), the flag module thus obtained only depends on the
class $gB\in G/B$.

Conversely, assume that $F$ is a flag module of $\O(\ol{G/U})$, and choose a
graded basis $\{e_\la\mid\la\in\Pp\}$ of $F$. For each $\la,\mu\in\Pp$,
let $v_\la^\mu\in V^\mu$ be defined by
\[
e_\la.f=\<f,v_\la^\mu\>\,e_{\la+\mu}\qquad\text{for all $f\in V_\mu$.}
\]
Since the algebra $\O(\ol{G/U})$ is commutative, we have
$(e_0.f).f'=(e_0.f').f$ for every $f\in V_\la$, $f'\in V_\mu$, hence
\begin{equation}
\label{symm0}
v_0^\la\o v_\la^\mu=v_\mu^\la\o v_0^\mu.
\end{equation}
It follows that $v_\la^\mu$ is a multiple of $v_0^\mu=:v^\mu$, say
$v_\la^\mu=a_\la v^\mu$. Inserting back into \eqref{symm0} yields
$a_\la=a_\mu$, for all $\la,\mu\in\Pp$. Since $a_0=1$, we get $a_\la=1$ for
all $\la\in\Pp$. Therefore,
\[
e_0.f=\<f,v^\mu\>\,e_\la\qquad\text{for all $f\in V_\la$.}
\]
The collection $\{v^\la\mid\la\in\Pp\}$ defines a linear form $v$ on
$\O(\ol{G/U})$. Furthermore, we have $e_0.(ff')=(e_0.f).f'$ for all $f\in
V_\la$, $f'\in V_\mu$, so $\<ff',v^{\la+\mu}\>=\<f,v^\la\>\<f',v^\mu\>$,
which shows that the linear form $v$ is a character on $\O(\ol{G/U})$,
corresponding to a point $x$ of the affine variety $\ol{G/U}$. Moreover,
since $F$ is cyclic, each $v^\la$ must be nonzero, so $x$ actually lies in
$G/U$, say $x=gU$. This yields an element $gB\in G/B$.

It is clear that these two constructions are inverse to each other.
\end{proof}
We will now discuss a possible picture of this kind in the noncommutative
situation: if $A$ is a quantum $G$, we would like to parametrize the
isomorphism classes of flag modules over the shape algebra $M_A$ by the
(closed) points of some scheme $E$.

Moreover, given a flag module $F$ and a weight $\la\in\Pp$, consider
the \emph{shifted} flag module $F[\la]$, defined as the $\Pp$-graded
module such that $F[\la]_\mu:=F_{\la+\mu}$. We would then like that, for
each $\la$, the shift operation $F\mapsto F[\la]$ corresponds to an
automorphism of schemes $\sig_\la:E\to E$.

To achieve this, let us encode the structure of a flag module more
geometrically, as follows. If $F$ is a flag module with basis
$\{e_\la\mid\la\in\Pp\}$, then for each $\la\in\Pp$ and each $1\le
i\le\ell$, let $v_\la^i\in V^i$ be defined by
\begin{equation}
\label{defmodstr}
e_\la.f=\<f,v_\la^i\>\,e_{\la+\om_i}\qquad\text{for all $f\in V_i$.}
\end{equation}
Replacing $F$ by an isomorphic flag module (i.e.\ rescaling the
$e_\la$) only multiplies each $v_\la^i$ by a scalar, so let
$p_\la^i$ be the corresponding point in $\P(V^i)$. To
simplify notation, let us write
\[
\MP:=\P(V^1)\times\dots\times\P(V^\ell)
\]
and denote by $\pr^i:\MP\to\P(V^i)$ the natural projection. For any
point $p\in\MP$, we use the shorthand notation $p^i:=\pr^i(p)$. Thus,
to an isomorphism class of flag modules, we associate a collection of
points $\{p_\la\mid\la\in\Pp\}$ in $\MP$.

{}From now on, we assume that Question~\ref{quadratic} has a positive
answer. The quadratic relations (I) and (II) in $M_A$ (see the end of
Section~\ref{shapesection}) impose some conditions on this collection of
points, which we now analyse.

For relations of type (I), identify $\P(V^i)\times\P(V^j)$ with its
image in $\P(V^i\o V^j)$ under the Segre embedding. Relations
(I)$_{ij}$ may be viewed as equations defining a subscheme
$\Gamma^{ij}$ of $\P(V^i)\times\P(V^j)$. We then have
\begin{equation}
\label{pinGamma}
(p_\la^i,p_{\la+\om_i}^j)\in\Gamma^{ij}
\end{equation}
for all $\la\in\Pp$ and all $1\le i,j\le\ell$.

Similarly, for relations of type (II), we consider the map
$\P(R^{ji}):\P(V^j\o V^i)\to\P(V^i\o V^j)$, where $R^{ji}$
denotes the transpose of $R_{ij}$. Then we must have
\begin{equation}
\label{pRij}
(p_\la^i,p_{\la+\om_i}^j)
=\P(R^{ji})(p_\la^j,p_{\la+\om_j}^i)
\end{equation}
(again identifying $\P(V^i)\times\P(V^j)$ with its image under the
Segre embedding).

Gluing together conditions~(\ref{pinGamma}) and~(\ref{pRij}) for all
$i,j$, we are
led to consider the subscheme $\Gamma\subset(\MP)^{\ell+1}$ of all
$(\ell+1)$-tuples $(p_0,p_1,\dots,p_\ell)$ satisfying
\begin{gather*}
(p_0^i,p_i^j)\in\Gamma^{ij},
\\
(p_0^i,p_i^j)=\P(R^{ji})(p_0^j,p_j^i)
\end{gather*}
for all $1\le i,j\le\ell$. We may now rephrase conditions (\ref{pinGamma})
and (\ref{pRij}) by saying that the collection $\{p_\la\mid\la\in\Pp\}$
satisfies
\begin{equation}
\label{Gammacond}
(p_\la,p_{\la+\om_1},\dots,p_{\la+\om_\ell})\in\Gamma
\qquad\text{for all $\la\in\Pp$.}
\end{equation}
\begin{prop}
\label{modulefamily}
Assume that $M_A$ is quadratic (as an $\N$-graded algebra). Then there
is a one-to-one correspondence between isomorphism classes of flag modules
over $M_A$ and families $\{p_\la\mid\la\in\Pp\}$ of points of $\MP$
satisfying (\ref{Gammacond}).
\end{prop}
\begin{proof}
It remains to show that the above construction can be reversed, so assume
that $\{p_\la\mid\la\in\Pp\}$ is a collection of points in $\MP$
satisfying (\ref{Gammacond}). Choose a (nonzero)
representative $v_\la^i\in V^i$ for each $p_\la^i$, and endow
a vector space $\bigoplus_{\la\in\Pp}\C e_\la$
with the flag module structure defined by the rule (\ref{defmodstr}). By
(\ref{pinGamma}), this rule is compatible with relations of type (I) in
$M_A$. By (\ref{pRij}), it is also compatible with relations of type (II),
\emph{provided} that, for each $\la\in\Pp$ and each $1\le i,j\le\ell$, we
suitably rescale one of $v_\la^i$, $v_{\la+\om_i}^j$, $v_\la^j$,
$v_{\la+\om_j}^i$. Proceeding by induction over the height $h(\la)$, we
may perform this rescaling in a consistent way.

It is clear that the two constructions are inverse to each other.
\end{proof}
\begin{rem}
Rescaling the $m_{\la\mu}$ only multiplies the $R_{ij}$ by scalars.
Therefore, the scheme $\Gamma$ does not depend on the normalizations of the
multiplication in $M_A$, but only on $A$ itself.
\end{rem}

The following Question is inspired by the description given in the
Introduction of~\cite{ATV1}.
\begin{quest}
\label{graphexists}
Do there exist a subscheme $E$ of $\MP$ and $\ell$ pairwise commuting
automorphisms $\sig_1,\dots,\sig_\ell:E\to E$ such that the scheme $\Gamma$
is given by
\begin{equation}
\label{Gammagraph}
\Gamma=\{(p,\sig_1(p),\dots,\sig_\ell(p))\mid p\in E\}?
\end{equation}
\end{quest}
A positive answer to this Question would fulfill the aim of parametrizing
flag modules, as expressed at the beginning of this Section. Indeed, assume
that $E$ and  $\sig_1,\dots,\sig_\ell$ as in Question~\ref{graphexists} do
exist. For each weight $\la=\sum a_i\om_i$, define
$\sig_\la:=\sig_1^{a_1}\dots\sig_\ell^{a_\ell}$; since the $\sig_i$ commute,
we have $\sig_{\la+\mu}=\sig_\la\sig_\mu$. Then for every family
$\{p_\la\mid\la\in\Pp\}$ of points in $\MP$ satisfying (\ref{Gammacond}), the
realization (\ref{Gammagraph}) shows that $p_\la=\sig_\la(p_0)$ for all
$\la\in\Pp$, with $p_0\in E$. Conversely, for every $p\in E$, the family
$\{\sig_\la(p)\mid\la\in\Pp\}$ satisfies (\ref{Gammacond}) and thus defines
an isomorphism class of flag modules by Proposition~\ref{modulefamily}.
Therefore, if Question~\ref{graphexists} had a positive answer, flag modules
(up to isomorphism) would be parametrized by the points of $E$, with
$\sig_\la$ corresponding to the shift operation $F\mapsto F[\la]$.

Finally, for future reference, we define, for each $1\le i\le\ell$, the line
bundle $\L_i$ over $E$ as the pullback of $\O_{\P(V^i)}(1)$ along $\pr^i$
(restricted to $E$), and we call the tuple
\[
\Triple(M_A):=(E,\sig_1,\dots,\sig_\ell,\L_1,\dots,\L_\ell)
\]
the \emph{flag tuple} associated to $A$.

\section{Braided tuples and reconstruction of shape algebras}
\label{Chansection}

Chan~\cite{Ch} has given a construction in the opposite direction, starting
from a scheme $E$, automorphisms $\sig_1,\dots,\sig_\ell$ of $E$, and line
bundles $\L_1,\dots,\L_\ell$ over $E$ (satisfying some compatibility
conditions; see Definition~\ref{tripledef}), and building a $\Pp$-graded
algebra from these data. Let us briefly recall his construction. To improve
legibility, we will write $\L^\sig$ for the pullback of a line bundle $\L$
along a map $\sig$.
\begin{df}
\label{tripledef}
We call a tuple $T=(E,\sigma_1,\dots,\sigma_\ell,\L_1,\dots,\L_\ell)$ as above
a \emph{braided tuple} if
\begin{enumerate}
\item the $\sigma_i$ pairwise commute,
\item for every $i>j$, there exists an equivalence
$R_{ij}:\L_i\o\L_j^{\sigma_i}\simto\L_j\o\L_i^{\sigma_j}$
of line bundles such that the braid relation
\begin{equation}
\label{Lijkbraid}
(R_{jk}\o\id)(\id\o R_{ik}^{\sig_j})(R_{ij}\o\id)
=(\id\o R_{ij}^{\sig_k})(R_{ik}\o\id)(\id\o R_{jk}^{\sig_i})
\end{equation}
holds for every $i>j>k$ (both sides being equivalences
$\L_i\o\L_j^{\sig_i}\o\L_k^{\sig_i \sig_j}\simto 
\L_k\o\L_j^{\sig_k}\o\L_i^{\sig_k \sig_j}$).
\end{enumerate}
\end{df}
Note that if we set $R_{ii}:=\id$ for all $i$ and $R_{ji}:=R_{ij}\inv$ for
all $i>j$, then \eqref{Lijkbraid} becomes true for all $i,j,k$.

If $\la\in\Pp$ decomposes as $\la=\sum a_i\om_i$, then define
$\sig_\la:=\sig_1^{a_1}\dots\sig_\ell^{a_\ell}$, as before (so
$\sig_{\la+\mu}=\sig_\la\sig_\mu$). Define a line bundle $\L_\la$ over
$E$ by the following inductive rules (with $\L_0$ the trivial bundle):
\begin{equation}
\label{Lrec}
\begin{gathered}
\L_{\om_i}=\L_i,\qquad 1\le i\le\ell,
\\
\L_{\la+\mu}=\L_\la\o\L_\mu^{\sig_\la}.
\end{gathered}
\end{equation}
As is shown in~\cite{Ch}, this procedure is, thanks to \eqref{Lijkbraid},
well defined up to unique equivalences of line bundles built from the
$R_{ij}$ (cf.\ also the proof of Proposition~\ref{assocbraid}). Now define
the product of two sections $\al\in \HZ(E,\L_\la)$ and $\be\in \HZ(E,\L_\mu)$
by
\begin{equation}
\label{Chanproduct}
\al\be:=\al\o\be^{\sig_\la}\in \HZ(E,\L_{\la+\mu}).
\end{equation}
\begin{theo}[Chan~\cite{Ch}]
The product rule~\eqref{Chanproduct} turns the direct sum
\[
B(T):=\bigoplus_{\la\in\Pp}\HZ(E,\L_\la)
\]
into an associative $\Pp$-graded algebra.
\end{theo}
The algebra $B(T)$ is not quadratic in general, so we consider its quadratic
cover
\[
M(T):=B(T)^{(2)}.
\]
(If $B$ is any $\N$-graded algebra, we define its \emph{quadratic cover}
$B^{(2)}$ as follows: consider the canonical homomorphism $\Tens(B_1)\to B$
and its kernel $J=\bigoplus_{k\ge2}J_k$, then set
$B^{(2)}:=\Tens(B_1)/(J_2)$. Here we view $B(T)$ as an $\N$-graded algebra
via the height function $h(\la)$.)

The quadratic algebra $M(T)$ may also be described more directly in terms of
the braided tuple $T$, as follows. For each $1\le i\le\ell$, set
$V_i:=\HZ(E,\L_i)$, denote by $V^i$ the dual of $V_i$, and consider the map
$\Pl^i:E\to\P(V^i)$ corresponding to the line bundle $\L_i$. For every $1\le
i,j\le\ell$, the map $\Pl^i\boxtimes(\Pl^j)^{\sig_i}$ corresponding to the line
bundle $\L_i\o\L_j^{\sig_i}$ is then given by the composite
\begin{equation}
\label{LiLjsi}
E\xto{\diag}E\times E\xto{\id\times\sig_i}
E\times E
\xto{\Pl^i\times\Pl^j}
\P(V^i)\times\P(V^j)
\xto{\Segre}\P(V^i\o V^j).
\end{equation}
Denote by $\Gamma^{ij}$ the image of this map and by $K_{ij}\subset
V_i\o V_j$ the subspace of linear forms on $V^i\o V^j$ vanishing on
$\Gamma^{ij}$.

For every $1\le i,j\le\ell$, Definition~\ref{tripledef} implies that there
exists a linear isomorphism $R^{ji}:V^j\o V^i\to V^i\o V^j$ such that the
following diagram commutes:
\begin{equation}
\label{crossdiag}
\xymatrix@R=3mm@C=1cm{
&E\times E\ar[r]^{\id\times\sig_i}
&E\times E\ar[r]^(.4){\Pl^i\times\Pl^j}
&\P(V^i)\times\P(V^j)\ar[r]^{\Segre}
&\P(V^i\o V^j)
\\
E\ar[dr]_(.4){\diag}\ar[ur]^(.4){\diag}
\\
&E\times E\ar[r]_{\id\times\sig_j}
&E\times E\ar[r]_(.4){\Pl^j\times\Pl^i}
&\P(V^j)\times\P(V^i)\ar[r]_{\Segre}
&\P(V^j\o V^i).\ar[uu]_{\P(R^{ji})}
}
\end{equation} 
Let $R_{ij}:V_i\o V_j\to V_j\o V_i$ be the transpose of $R^{ji}$. It is clear
that modulo $K_{ij}$ and $K_{ji}$, the map $R_{ij}$ is unique up to a scalar.

The algebra $M(T)$ is then generated by $V_1\oplus\dots\oplus V_\ell$,
with relations given by (I)$_{ij}$ and~(II)$_{ij}$ for all $1\le
i,j\le\ell$ (see the end of Section~\ref{shapesection};
Remark~\ref{superfluous} still applies).
\begin{quest}
What is the kernel of the canonical morphism $M(T)\to B(T)$?
\end{quest}
Having constructed the algebra $M(T)$ from a braided tuple $T$, we
may formulate a converse to Question~\ref{graphexists}:
\begin{quest}
\label{tripleexists}
Assume that $A$ is a quantum $G$ such that the shape algebra $M_A$ is
quadratic. Does there exist a braided tuple $T$ such that $M_A=M(T)$?
\end{quest}
This Question is \emph{a priori} weaker than Question~\ref{graphexists}, for
the following reason. If $M_A$ is quadratic and does admit a flag tuple $T$
as in Question~\ref{graphexists}, then the reconstructed algebra $M(T)$ is
canonically isomorphic to $M_A$. However, we might also have $M(T')=M_A$ for
some subtuple $T'$ of $T$ (i.e.\ a subscheme $E'$ of $E$ stabilized by each
$\sig_i$, with $\sig'_i$ and $\L'_i$ the obvious restrictions).
\begin{prob}
\label{intrinsic}
Given a simple complex group $G$, characterize the flag tuples of all
quantum $G$'s intrinsically (i.e.\ as braided tuples).
\end{prob}
For $G=\SL(2)$, this is elementary: $E$ must be the projective line $\P^1$,
$\sig$ can be an arbitrary automorphism of infinite order, and
$\L=\O_{\P^1}(1)$. The three possible forms of $\sig$ correspond to three
different quantum $\SL(2)$'s, namely, $\O(\SL(2))$ (when $\sig=\id$), the
standard Drinfel$'$d-Jimbo quantum $\SL(2)$ for $q$ not a root of unity (when
$\sig$ has two fixed points), and the Jordanian quantum $\SL(2)$~\cite{DMMZ}
(when $\sig$ has one fixed point). These are known~\cite{Wo} to be the only
quantum $\SL(2)$'s (in the sense of Definition~\ref{Gqdef}). The associated
shape algebras are $\C\<x,y\>/(xy-yx)$, $\C\<x,y\>/(xy-q\,yx)$, and
$\C\<x,y\>/(xy-yx-y^2)$, respectively.
\begin{prob}
Reconstruct not only a shape algebra, but a quantum $G$ itself from a
braided tuple satisfying the conditions found in
Problem~\ref{intrinsic}.
\end{prob}

\part{A conjectural flag tuple for the standard Drinfel$\mathbf'$d-Jimbo
quantum groups}
\label{djpart}

In this Part, we describe ingredients for a potential braided tuple, and we
conjecture that these geometric data provide positive answers to Questions
\ref{tripleexists}~and \ref{graphexists} for the standard quantum groups of
Drinfel$'$d and Jimbo. (The conjecture concerning Question~\ref{tripleexists}
will be proved for $\SL(n)$ in Part~\ref{djslnpart}.)

Again, $G$ will denote a simply-connected simple complex group.

\section{Recollections on $\Uqg$ and $\OqG$}
\label{DJsection}

Let $\g$ be the Lie algebra of $G$. Drinfel$'$d~\cite{Dr} and Jimbo~\cite{Ji}
have defined (independently) a Hopf algebra $\Uqg$ that depends on a
parameter $q\in\C^*$ and that is a ``quantum analogue'' of the universal
enveloping algebra $\U(\g)$ (in the sense that its comultiplication is no
longer cocommutative). When $q$ is not a root of unity, finite-dimensional
$\Uqg$-mod\-ules have been studied (independently) by Lusztig~\cite{Lu} and
by Rosso~\cite{Ro}: in particular, discarding unwanted nontrivial
one-dimensional modules, there still exists a family $\{V_\la\mid\la\in\Pp\}$
of $\U_q(\g)$-mod\-ules satisfying conditions (a) and~(c) of
Definition~\ref{Gqdef} (condition (c) follows e.g.\ from Theorem~4.12(b)
of~\cite{Lu}).

Therefore, if $\OqG$ denotes the subspace of $\Uqg^*$ spanned by the matrix
coefficients of the modules $V_\la$, then $\OqG$ (for $q$ not a root of
unity) is a quantum $G$ in the sense of Definition~\ref{Gqdef}. (If $G$ is
not assumed to be simply connected, then $\OqG$ may still be defined in this
way, provided $\Pp$ is replaced by the appropriate submonoid.) We call $\OqG$
the \emph{standard} quantum $G$. When $G=\SL(n)$, $\SO(n)$, or $\Sp(n)$, a 
presentation of $\OqG$ by generators and relations has been given by Faddeev,
Reshetikhin, and Takhtajan~\cite{FRT}.

\section{Recollections from~\cite{PPinGB}}
\label{sectPPinGB}

Choose a Borel subgroup $B\subset G$ and a maximal torus $T\subset B$, and
let $W:=\Norm_G(T)/T$ be the associated Weyl group. Denote by $\Phi$ and
$\Phi^+$ the root system and the set of positive roots, respectively. To each
$\al\in\Phi$ are associated a reflection $s_\al\in W$, a root group $U_\al$,
and a copy $L_\al=\<U_\al,U_{-\al}\>$ of $\operatorname{(P)SL}(2)$ in $G$.

Recall the following construction from~\cite{PPinGB}: an \emph{orthocell}
(of rank $d$) is a left coset in $W$ of the form
\[
C=C(w;\al_1,\dots,\al_d):=w\<s_{\al_1},\dots,s_{\al_d}\>,
\]
where $w\in W$ and $\al_1,\dots,\al_d$ are positive and pairwise orthogonal
roots.

\emph{Warning:} the $\al_k$ are \emph{not} assumed to be \emph{strongly}
orthogonal, i.e.\ the sum of two of them may well be a root.

By orthogonality, the reflections $s_{\al_1},\dots,s_{\al_d}$ pairwise
commute. Therefore, the following notation makes sense, and we will use it
frequently:
\[
s_L:=\prod_{k\in L}s_{\al_k},\qquad L\subset\{1,\dots,d\}.
\]
(Note that the elements of $C$ are those of the form $ws_L$.)

Reordering the sequence $\al_1,\dots,\al_d$ if necessary, assume that it is
\emph{nonincreasing}, in the sense that $\al_k\not<\al_{k'}$ for all $k<k'$;
then define
\[
E(C):=\{\dot wg_1\dots g_dB\mid g_k\in L_{\al_k}\ \forall k\}\subset G/B,
\]
where $\dot w\in\Norm_G(T)$ is some representative of $w$. In~\cite{PPinGB},
we show that $E(C)$ only depends on $C$ as a coset (and not on the choice of
$w$ in $C$, nor of its representative $\dot w$, nor on the chosen
nonincreasing ordering of the $\al_k$). Furthermore, we show that $E(C)$ is a
$T$-stable subvariety of $G/B$, isomorphic to the product
$\P^1\times\dots\times\P^1$ of $d$ projective lines.

\begin{rem}
\label{rightrem}
Orthocells may also be defined in terms of \emph{right} cosets: if we set
\[
C(\al_1,\dots,\al_d;w):=\<s_{\al_1},\dots,s_{\al_d}\>w,
\]
then $C(\al_1,\dots,\al_d;w)=C(w;w\inv\al_1,\dots,w\inv\al_d)$. Moreover,
recall (see~\cite{Sp}, end of~\S9.2.1) that for each $w\in W$ and each root
$\al$, we have $wU_\al w\inv=U_{w\al}$, and hence $wL_\al w\inv=L_{w\al}$. It
follows that for $C=C(\al_1,\dots,\al_d;w)$, we have
\[
E(C)=\{g_1\dots g_dwB\mid g_k\in L_{\al_k}\ \forall k\}.
\]
\end{rem}

\section{Monogressive orthocells and the variety $E^\DrJ$}
\label{EDJG}

Denote by $<$ the Bruhat order on $W$ and by $\lessdot$ the associated cover
relation (i.e.\ $w\lessdot w'$ if $w<w'$ and if no element of $W$ lies between
$w$ and $w'$). Denote also by $\ell(w)$ the length of an element $w\in W$.
Recall the following combinatorial characterization (see e.g.\ Sections 5.9
and~5.11 of~\cite{Hu}):
\[
w\lessdot w'\iff\ell(w')=\ell(w)+1\text{ and }w'=ws\text{ for some reflection
$s$.}
\]
Assume that $w\in C$ has been chosen of minimal length.
\begin{df}
An orthocell $C=C(w;\al_1,\dots,\al_d)$ will be called \emph{monogressive} if
\[
ws_L\lessdot ws_Ls_{\al_k}\qquad\forall L\subset\{1,\dots,d\},
\quad\forall k\not\in L,
\]
or, equivalently, if
\[
\ell(ws_L)=\ell(w)+|L|\qquad\forall L\subset\{1,\dots,d\}.
\]
\end{df}
We then define the variety $E^\DrJ\subset G/B$ by
\[
E^\DrJ:=\bigcup_{\text{$C$ monogressive}} E(C).
\]

\section{The automorphisms $\sig_1,\dots,\sig_\ell$}
\label{sigDJG}

Let $\be_1,\dots,\be_\ell$ be the simple roots. Then the morphism
\[
T\to(\C^*)^\ell:t\mapsto(\be_1(t),\dots,\be_\ell(t))
\]
is surjective, so we may choose, for each $1\le i\le\ell$, an element $t_i\in
T$ such that
\[
\be_j(t_i)=q^{-(\om_i|\be_j)}=
\begin{cases}
q^{-(\be_j|\be_j)/2}&\text{if $j=i$},
\\
1&\text{if $j\ne i$.}
\end{cases}
\]
By (multiplicative) linearity, it then follows that
$\al(t_i)=q^{-(\om_i|\al)}$ for every root $\al\in\Phi$.

Now let $C=C(w;\al_1,\dots,\al_d)$ be a monogressive orthocell. Since $E(C)$ is
$T$-stable in $G/B$, the automorphism
\[
\sig_{i,C}:E(C)\to E(C):gB\mapsto wt_i w\inv gB
\]
is well defined, and it is independent of the choice of $t_i$ because the
kernel of the above morphism $T\to(\C^*)^\ell$ is equal to the centre of $G$
(see e.g.~\cite[Proposition~8.1.1]{Sp}).
\begin{prop}
\label{glue}
For each $i$, the automorphisms $\sig_{i,C}$ glue together to form a well
defined automorphism $\sig_i$ of $E^\DrJ$.
\end{prop}
\noindent We defer the proof to Appendix~\ref{glueproof}.

\section{The line bundles $\L_1,\dots,\L_\ell$}
\label{LDJG}

Recall that for each $\la\in\Pp$, the highest weight point in $\P(V^\la)$ is
fixed by $B$, hence we get a well defined \emph{Pl\"ucker map}
\[
\Pl^\la:G/B\to\P(V^\la).
\]
Let $\om_1,\dots,\om_\ell$ be the fundamental weights and write
$\Pl^i:=\Pl^{\om_i}$ for each $i$. We then define the a line
bundle $\L_i$ as the pullback of $\O_{\P(V^i)}(1)$ along $\Pl^i$, restricted
to $E^\DrJ$.

\emph{Warning.} We may not define $\L_\la$ to be the pullback of
$\O_{\P(V^\la)}(1)$ for all $\la\in\Pp$: this would cause a conflict with the
recursion rule \eqref{Lrec}.

\section{Main conjectures and result}

\begin{conj}[Positive answer to Question~\ref{tripleexists}]
\label{mainconj}
Assume that $q\in\C^*$ is not a root of unity. The tuple
$T^\DrJ=(E^\DrJ,\sig_1,\dots,\sig_\ell,\L_1,\dots,\L_\ell)$ defined in
Sections \ref{EDJG}--\ref{LDJG} is a braided tuple (see
Definition~\ref{tripledef}), and the associated quadratic algebra $M(T^\DrJ)$
is the shape algebra of the standard quantum group $\OqG$.
\end{conj}
\begin{conj}[Positive answer to Question~\ref{graphexists}]
\label{secconj}
Moreover, the same tuple $T^\DrJ$ is the flag tuple associated to $\OqG$
(i.e.\ $E^\DrJ$ parametrizes \emph{all} flag modules of the shape algebra of
$\OqG$).
\end{conj}
\begin{theo}
Conjecture~\ref{mainconj} is true for $G=\SL(n)$.
\end{theo}

\part{The standard quantum $\SL(n)$}
\label{djslnpart}

In this Part, we will describe the objects of Sections
\ref{sectPPinGB}--\ref{sigDJG} more explicitely when $G=\SL(n)$, and we prove
Conjecture~\ref{mainconj} in that case.

\section{The varieties $E(C)$}

From now on, it will be more convenient to view orthocells as \emph{right}
cosets (see Remark~\ref{rightrem}).

Let us first recall the usual realization of the flag variety $\SL(n)/B$, of
the Pl\"ucker maps $\Pl^i$, and of the subgroups $L_\al$.

We let $B\subset\SL(n)$ be the subgroup of all upper triangular matrices,
i.e.\ the stabilizer of the flag
\[
\C e_1\subset\C e_1\oplus\C e_2\subset\dots\subset\C e_1\oplus\dots\oplus\C
e_{n-1},
\]
where $e_1,\dots,e_n$ denotes the canonical basis of $\C^n$. This identifies
$\SL(n)/B$ with the set of all (full) flags in $\C^n$ (or in
$\P^{n-1}:=\P(\C^n)$).

We also let $T\subset B$ be the subgroup of all diagonal matrices: the Weyl
group $W$ then identifies with the symmetric group $S_n$, and the reflections
correspond exactly to the transpositions.

For each $1\le i\le n-1$, recall that the fundamental representation
$V^i:=V^{\om_i}$ is given by the exterior power $\Lambda^i\C^n$, and that the
map $\Pl^i:\SL(n)/B\to\P(\Lambda^i\C^n)$ may be described as
follows: given a flag $F\in\SL(n)/B$, choose a basis $f_1,\dots,f_i$ of its
component $F_i$ of dimension $i$, then send $F$ to the point
$f_1\wedge\dots\wedge f_i\in\P(\Lambda^i\C^n)$ (which is independent of
the choice of the basis). Moreover, the elements
\[
e_w^i:=e_{w(1)}\wedge\dots\wedge e_{w(i)},\qquad w\in S_n,
\]
form a basis of $\Lambda^i\C^n$ (up to obvious redundancies). 

Let $\al\in\Phi^+$ and write $s_\al=(a\,b)$, $1\le a<b\le n$. Then the
subgroup $L_\al\subset\SL(n)$ is the group $\SL(2)$ acting naturally on $\C
e_a\oplus\C e_b$ and trivially on all other $e_c$. Clearly, if $s_\al,s_\be$
commute (i.e.\ if $\al,\be$ are orthogonal), then so do $L_\al$ and $L_\be$.
(This is not true for arbitrary $G$.)

Now fix an orthocell $C=C(\al_1,\dots,\al_d;w)$ and let us describe the
variety $E(C)$, or rather, its images under the maps $\Pl^i$, $1\le i\le
n-1$.
\begin{rem}
\label{equiveff}
For each $1\le k\le d$, the following conditions are equivalent:
\begin{itemize}
\item $s_{\al_k}w\om_i=w\om_i$,
\item the transposition $s_{\al_k}$ leaves the set $\{w(1),\dots,w(i)\}$
invariant,
\item $e_{s_{\al_k}w}^i=\pm e_w^i$.
\end{itemize}
\end{rem}
Number the $\al_k$ in such a way that for some $1\le a\le d$, the above
conditions hold for $1\le k\le a$ and do not hold for $a+1\le k\le d$. For
each $k$, write $s_{\al_k}=(a_k\,b_k)$, $a_k<b_k$, and pick an
element $g_k\in L_{\al_k}$ acting as $\spmat{x_k&*\\y_k&*}$ on $\C
e_{a_k}+\C e_{b_k}$ (and trivially on the other $e_c$). For any subset
$L\subset\{1,\dots,a\}$, write
\[
\bar L:=\{1,\dots,a\}\setminus L,
\qquad
x_L:=\prod_{k\in L}x_k,
\qquad
y_L:=\prod_{k\in L}y_k.
\]
The above description of the map $\Pl^i$ and of the subgroups
$L_{\al_k}$ now imply that
\begin{equation}
\label{Flexpr}
\Pl^i(g_1\dots g_dwB)=\sum_{L\subset\{1,\dots,a\}} x_{\bar L}y_L\,
e_{s_Lw}^i\in\P(\Lambda^i\C^n).
\end{equation}
\begin{rem}
\label{homcoord}
The variety $E(C)$ being a product of $d$ projective lines, we may view
$(x_1:y_1),\dots,(x_d:y_d)$ as homogeneous coordinates on these lines.
\end{rem}
A more geometric description of the varieties $E(C)$ (not needed here) can be
found in~\cite[Example~5.1]{PPinGB}.

\section{Monogressivity}

Let us first recall a more explicit description of the Bruhat cover relation
in $S_n$. Write a permutation $w\in S_n$ as an array $[w(1)\dots w(n)]$, and
write e.g.\ $w=[\cdot a\cdot b\cdot c\cdot]$ to signify that in the array
$w$, $a$ appears to the left of $b$ and $b$ appears to the left of $c$.

If $s\in S_n$ is a transposition, say $s=(a\,b)$ with $a<b$, then 
$w\lessdot sw$ if and only if (i)~$w=[\cdot a\cdot b\cdot]$ and (ii)~whenever
$w=[\cdot a\cdot c\cdot b\cdot]$, $c$ is outside of the (numerical) interval
$[a,b]$. For example, if $n=7$ and $s=(4\,6)$, then we have
$[3472651]\lessdot[3672451]$, but $[7415623]{\,\not\!\!\lessdot\,}[7615423]$
because the subarray $[4156]$ contains $5$.

Now let $C=C(\al_1,\dots,\al_d;w)$ be an orthocell, and write again
$s_{\al_k}=(a_k\,b_k)$, $a_k<b_k$, for all $k$. The above description of the
Bruhat cover relation shows the following
\begin{crit}
\label{monocrit}
With the above notation, the orthocell $C$ is monogressive if and only if the
following conditions hold for all $k$:
\begin{itemize}
\item $w=[\cdot a_k\cdot b_k\cdot]$,
\item whenever $w=[\cdot a_k\cdot c\cdot b_k\cdot]$, $c$ is outside of the
interval $[a_k,b_k]$,
\item whenever $w=[\cdot a_k\cdot a_{k'}\cdot b_k\cdot]$ for some $k'\ne k$,
both $a_{k'}$ and $b_{k'}$ are outside of the interval $[a_k,b_k]$, and
similarly whenever $w=[\cdot a_k\cdot b_{k'}\cdot b_k\cdot]$.
\end{itemize}
\end{crit}
\begin{ex}[$n=4$]
There are fifty-eight monogressive orthocells of rank~$1$, viz.\ those of
one of the following forms:
\begin{itemize}
\item $\{[ijkl],[jikl]\}$, $\{[kijl],[kjil]\}$, or $\{[klij],[klji]\}$ ($i<j$),
\item $\{[ikjl],[jkil]\}$ or $\{[likj],[ljki]\}$ ($i<j$; $k\not\in[i,j]$),
\item $\{[iklj],[jkli]\}$ ($i<j$; $k,l\not\in[i,j]$).
\end{itemize}
There are eleven monogressive orthocells $C(\al_1,\al_2;w)$ of rank~$2$,
given by
\begin{itemize}
\item $s_{\al_1}=(1\,2)$, $s_{\al_2}=(3\,4)$, $w=[1234]$,
$[3412]$, $[1324]$, $[3142]$, $[1342]$, or $[3124]$;
\item $s_{\al_1}=(1\,3)$, $s_{\al_2}=(2\,4)$, $w=[1324]$ or
$[2413]$;
\item $s_{\al_1}=(1\,4)$, $s_{\al_2}=(2\,3)$, $w=[1423]$,
$[2314]$, or $[2143]$.
\end{itemize}
(Pictures for the corresponding varieties $E(C)$ may be found
in~\cite[Example~5.1]{PPinGB}.)
\end{ex}
\begin{ex}[$n=3$]
\label{EforSL3}
View $\SL(3)/B$ as the set of flags $(p,l)$ in $\P^2=\P(\C^3)$.
Consider $e_1,e_2,e_3$ as points in $\P^2$ and let $e_{ab}\subset\P^2$ be the
line through $e_a$ and $e_b$. Then $E^\DrJ$ is the union of the following eight
curves in $\SL(3)/B$:
\begin{align*}
&\{(e_a,l)\mid e_a\in l\},\quad a=1,2,3,
\\
&\{(p,e_{ab})\mid p\in e_{ab}\},\quad ab=12,13,23,
\\
&\{(p,l)\mid e_{12}\ni p\in l\ni e_3\},
\\
&\{(p,l)\mid e_{23}\ni p\in l\ni e_1\}.
\end{align*}
See Figure~\ref{SL3picture}.
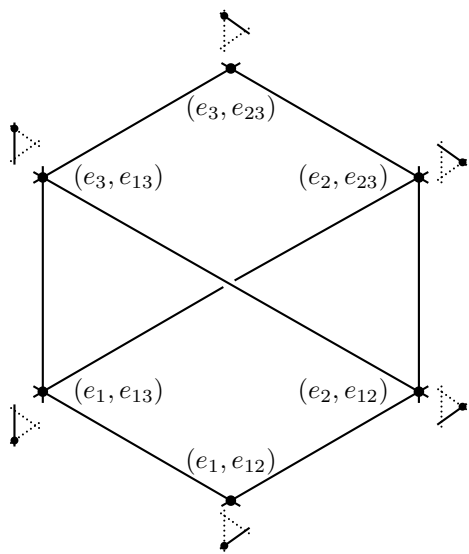
\begin{figure}
\begin{center}
\psset{unit=.5mm}
\begin{pspicture}(0,-10)(100,135)
\qline(47.5,-1.45)(102.5,30.45)
\qline(100,26)(100,89)
\qline(102.5,84.55)(47.5,116.45)
\qline(52.5,116.45)(-2.5,84.55)
\qline(0,89)(0,26)
\qline(-2.5,30.45)(52.5,-1.45)
\psline(-2.5,27.55)(48,56.84)
\psline(51,58.58)(102.5,87.45)
\psline(-2.5,87.45)(102.5,27.55)
\pscircle*(50,0){2pt}
\rput[b]{0}(50,7){$(e_1,e_{12})$}
\pscircle*(100,29){2pt}
\rput[r]{0}(92,29){$(e_2,e_{12})$}
\pscircle*(100,86){2pt}
\rput[r]{0}(92,86){$(e_2,e_{23})$}
\pscircle*(50,115){2pt}
\rput[t]{0}(50,108){$(e_3,e_{23})$}
\pscircle*(0,86){2pt}
\rput[l]{0}(8,86){$(e_3,e_{13})$}
\pscircle*(0,29){2pt}
\rput[l]{0}(8,29){$(e_1,e_{13})$}
\rput[t]{0}(51,-4){\psset{unit=.4mm}
\begin{pspicture}(0,0)(7,10)
\psline[linestyle=dotted,dotsep=1pt](0,-2)(0,12)
\psline[linestyle=dotted,dotsep=1pt](-1.4,11)(8.4,4)
\psline(8.4,6)(-1.4,-1)
\pscircle*(0,0){1.5pt}
\end{pspicture}}
\rput[tl]{0}(106,29){\psset{unit=.4mm}
\begin{pspicture}(0,0)(7,10)
\psline[linestyle=dotted,dotsep=1pt](0,-2)(0,12)
\psline[linestyle=dotted,dotsep=1pt](-1.4,11)(8.4,4)
\qline(8.4,6)(-1.4,-1)
\pscircle*(7,5){1.5pt}
\end{pspicture}}
\rput[tr]{0}(-2,24){\psset{unit=.4mm}
\begin{pspicture}(0,0)(7,10)
\psline[linestyle=dotted,dotsep=1pt](8.4,6)(-1.4,-1)
\psline[linestyle=dotted,dotsep=1pt](-1.4,11)(8.4,4)
\qline(0,-2)(0,12)
\pscircle*(0,0){1.5pt}
\end{pspicture}}
\rput[bl]{0}(106,86){\psset{unit=.4mm}
\begin{pspicture}(0,0)(7,10)
\psline[linestyle=dotted,dotsep=1pt](0,-2)(0,12)
\psline[linestyle=dotted,dotsep=1pt](8.4,6)(-1.4,-1)
\qline(-1.4,11)(8.4,4)
\pscircle*(7,5){1.5pt}
\end{pspicture}}
\rput[br]{0}(-2,91){\psset{unit=.4mm}
\begin{pspicture}(0,0)(7,10)
\psline[linestyle=dotted,dotsep=1pt](8.4,6)(-1.4,-1)
\psline[linestyle=dotted,dotsep=1pt](-1.4,11)(8.4,4)
\qline(0,-2)(0,12)
\pscircle*(0,10){1.5pt}
\end{pspicture}}
\rput[b]{0}(51,121){\psset{unit=.4mm}
\begin{pspicture}(0,0)(7,10)
\psline[linestyle=dotted,dotsep=1pt](0,-2)(0,12)
\psline[linestyle=dotted,dotsep=1pt](8.4,6)(-1.4,-1)
\qline(-1.4,11)(8.4,4)
\pscircle*(0,10){1.5pt}
\end{pspicture}}
\end{pspicture}
\end{center}
\caption{The subvariety $E^\DrJ$ in $\SL(3)/B$. Its eight irreducible
components intersect in six points. Next to each point is a 
small picture, viewing it as a flag in $\P^2$.
The ``missing'' diagonal corresponds to the orthocell
$\{[123],[321]\}$, which is not monogressive.}
\label{SL3picture}
\end{figure}
\end{ex}

\section{The automorphisms $\sig_i$}

Denote again  the simple roots by $\be_1,\dots,\be_{n-1}$. If
$t=\operatorname{diag}(x_1,\dots,x_n)\in T$ (with $\prod_j x_j=1$), then
recall that $\be_i(t)=x_ix_{i+1}\inv$. Therefore, $t_i\in T$ is equal, up to a
factor, to the matrix $\operatorname{diag}(1,\dots,1,q,\dots,q)$ ($i$
times $1$ and $n-i$ times $q$).

If $C=C(\al_1,\dots,\al_d;w)$, with $s_{\al_k}=(a_k\,b_k)$ as before, then
the action of the associated automorphism $\sigma_i:gB\mapsto wt_iw\inv gB$ on
$E(C)$ may be described more explicitly using the homogeneous coordinates
$(x_1:y_1),\dots,(x_d:y_d)$ of Remark~\ref{homcoord}:
\[
\sigma_i:(x_k:y_k)\mapsto
\begin{cases}
(x_k:qy_k)&\text{if $s_{\al_k}w\om_i\ne w\om_i$,}
\\
(x_k:y_k)&\text{if $s_{\al_k}w\om_i=w\om_i$.}
\end{cases}
\]
\begin{ex}[$n=3$]
On each of the eight components of $E^\DrJ$ (see Example~\ref{EforSL3}),
$\sig_1,\sig_2$ act as homotheties (viewing the two $T$-stable points on this
component as $0$ and $\infty$). The ratios for $\sig_1$ are, respectively,
$1,1,1,q,q,q,q,q$, and those for $\sig_2$ are $q,q,q,1,1,1,q,q$.
\end{ex}

\section{Proof of Conjecture~\ref{mainconj} for $G=\SL(n)$}

For each $1\le i\le n-1$, consider the vector space $V^i:=\Lambda^i\C^n$. On
one hand, $\SL(n)$ acts on it naturally, and the corresponding map
$\SL(n)/B\to\P(V^i)$ induces a line bundle $\L_i$ on $E^\DrJ\subset\SL(n)/B$
(see Section~\ref{LDJG}). On the other hand, we will make $\Uqsln$ act on
$V^i$ (see below, before Lemma~\ref{ijinv}), turning $V^i$ into the simple
$\Uqsln$-mod\-ule of highest weight $\om_i$.

Both the algebra $M(T^\DrJ)$ and the shape algebra $M^\DrJ$ thus become
quotients of the tensor algebra $T(V_1\oplus\dots\oplus V_{n-1})$ (where
$V_i:=(V^i)^*$). Note also that $M(T^\DrJ)$ is quadratic by definition, and
$M^\DrJ$ is quadratic by \cite{TT}. So we need to show that relations of
types (I)~and (II) (see end of Section~\ref{shapesection}) agree for both
algebras (and, of course, that the tuple $T^\DrJ$ is braided in the first
place).

We will break down the proof into several lemmas.

\begin{df}
Let $1\le i,j\le n-1$. An orthocell $C=C(\al_1,\dots,\al_d;w)$ will be called
\emph{$ij$-ef\-fect\-ive} if, for every $1\le k\le d$, we have both
$s_{\al_k}w\om_i\ne w\om_i$ and $s_{\al_k}w\om_j\ne w\om_j$. In this case, we
define the following element of $V^i\o V^j$:
\[
e_C^{ij}
:=\sum_{L\subset\{1,\dots,d\}} q^{|L|} e_{s_{\bar L}w}^i \o
e_{{s_L}w}^j,
\]
where, as before, $s_L:=\prod_{k\in L}s_{\al_k}$ and $\bar
L:=\{1,\dots,d\}\setminus L$.
\end{df}
We denote by $V^{ij}\subset V^i\o V^j$ the linear span of the image of the
map $\Pl^i\boxtimes(\Pl^j)^{\sig_i}$ (see~\eqref{LiLjsi}).
\begin{lem}
\label{spanned}
The subspace $V^{ij}$ is linearly spanned by the $e_C^{ij}$ for $C$
monogressive and $ij$-ef\-fect\-ive.
\end{lem}
\begin{proof}
First, let $C=C(\al_1,\dots,\al_d;w)$ be monogressive and $ij$-ef\-fect\-ive.
If a point $p\in E(C)$ has homogeneous coordinates
$(x_1:y_1),\dots,(x_d:y_d)$ (see Remark~\ref{homcoord}), then by
$ij$-ef\-fect\-ive\-ness, the coordinates of $\sig_i(p)$ are
$(x_1:qy_1),\dots,(x_d:qy_d)$. Using~\eqref{Flexpr}, we therefore see that
$\Pl^i\boxtimes(\Pl^j)^{\sig_i}$ sends $p$ to the following point in $\P(V^i\o
V^j)$:
\[
\sum_{L,M\subset\{1,\dots,d\}}q^{|M|}x_{\bar L}x_{\bar M}y_Ly_M\,
e_{s_L w}^i\o e_{s_M w}^j.
\]
Using the change of ``variables'' $I:=L\triangle M$ (symmetric difference),
$J:=L\cap M$, and $N:=M\setminus L=M\setminus J\subset I$, this expression
may be rewritten as
\begin{align*}
&\sum_{\substack{I,J\subset\{1,\dots,d\}\\I\cap J=\emptyset}}
q^{|J|} (x_{\ol{I\cup J}})^2 x_I y_I (y_J)^2
\left(\sum_{N\subset I} q^{|N|} e_{s_{\bar N}s_Jw}^i\o e_{s_Ns_Jw}^j
\right)
\\
&\qquad=\sum_{\substack{I,J\subset\{1,\dots,d\}\\I\cap J=\emptyset}}
q^{|J|} (x_{\ol{I\cup J}})^2 x_I y_I (y_J)^2\,
e_{C(\al_I;s_Jw)}^{ij},
\end{align*}
where $\al_I$ is shorthand for the set $\{\al_k\mid k\in I\}$. By induction
over $d$, we may assume that $e_{C'}^{ij}\in V^{ij}$ for all monogressive
$ij$-ef\-fect\-ive orthocells $C'$ of rank smaller than $d$ (the case $d=0$
being trivial). Since the above sum is in $V^{ij}$ by definition, the only
remaining term, namely $e_C^{ij}$, is in $V^{ij}$ as well.

We still need to show that the image of a point $p\in E(C)$ is in the span of
the $e_{C'}^{ij}$ (for $C'$ monogressive and $ij$-ef\-fect\-ive) even if $C$
is not $ij$-ef\-fect\-ive (but still monogressive). Reordering the $\al_k$ if
necessary, we may assume that, for some $1\le a\le b\le c\le d$, they satisfy
\begin{align*}
s_{\al_k}w\om_i&\ne w\om_i,&s_{\al_k}w\om_j&\ne w\om_j&&\text{if
$\phantom{a+{}}1\le k\le
a$;}
\\
s_{\al_k}w\om_i&\ne w\om_i,&s_{\al_k}w\om_j&=w\om_j&&\text{if $a+1\le k\le
b$;}
\\
s_{\al_k}w\om_i&=w\om_i,&s_{\al_k}w\om_j&\ne w\om_j&&\text{if $b+1\le k\le
c$;}
\\
s_{\al_k}w\om_i&=w\om_i,&s_{\al_k}w\om_j&=w\om_j&&\text{if $c+1\le k\le
d$.}
\end{align*}
Let again $(x_1:y_1),\dots,(x_d:y_d)$ be homogeneous coordinates for a point
$p\in E(C)$. This time, the coordinates for $\sig_i(p)$ are obtained by
multiplying $y_k$ by $q$ only for $1\le k\le b$. Furthermore, we have
\[
\Pl^i(p)=
\sum_{\substack{L\subset\{1,\dots,a\}\\L'\subset\{a+1,\dots,b\}}}
x_{\bar L} x_{\bar L'} y_L y_{L'}\,e_{s_Ls_{L'}w}^i,
\]
and a similar expression for $\Pl^j(p)$, with $\{a+1,\dots,b\}$ replaced by
$\{b+1,\dots,c\}$. A computation similar to the one above shows that
$\Pl^i\boxtimes(\Pl^j)^{\sig_i}$ now sends $p$ to
\begin{equation}
\label{imagep}
\sum_{\substack{I,J\subset\{1,\dots,a\}\\I\cap J=\emptyset\\
L'\subset\{a+1,\dots,b\}\\M'\subset\{b+1,\dots,c\}}}
q^{|J|} (x_{\ol{I\cup J}})^2 x_I x_{\bar L'} x_{\bar M'} y_I (y_J)^2 y_{L'}
y_{M'}\, e_{C(\al_I;s_Js_{L'}s_{M'}w)}^{ij}
\end{equation}
(where we have used the fact that $e_w^i=e_{s_{M'}w}^i$ and
$e_w^j=e_{s_{L'}w}^j$).
\end{proof}
\begin{lem}
\label{rightnumber}
Let $1\le i,j\le n-1$, with, say, $i<j$. Then
\[
\dim V^{ij}\le D_{n;i,j}:=\binom ni\binom nj-\binom n{i-1}\binom n{j+1}.
\]
\end{lem}
\begin{proof}
Consider a monogressive $ij$-ef\-fect\-ive cell $C=C(\al_1,\dots,\al_d;w)$.
For each $1\le k\le d$, write again $s_{\al_k}=(a_k\,b_k)$, $a_k<b_k$ (hence
$w=[\cdot a_k\cdot b_k\cdot]$ by monogressivity). Reorder the $\al_k$ in such
a way that $b_1<\dots<b_d$. By $ij$-ef\-fect\-ive\-ness, each $a_k$ must
appear in the subarray $[w(1)\dots w(i)]$, and each $b_k$ in the subarray
$[w(j+1)\dots w(n)]$. Now let $S_{ij}:=S_i\times S_{j-i}\times
S_{n-(i+j)}\subset S_n$, and note that for each $\pi\in S_{ij}$, replacing
$w$ by $w\pi$ in $C$ leaves $e_C^{ij}$ invariant up to a sign. Choosing $\pi$
appropriately, we may assume that $w$ takes the following form:
\[
[w(1)\dots w(i-d)\,
a_1\dots a_d\,
w(i+1)\dots w(j)\,
b_d\dots b_1\, 
w(j+d+1)\dots w(n)],
\]
with, say, the following orderings:
\[
w(1)>\dots>w(i-d),\qquad w(i+1)>\dots>w(j),\qquad w(j+d+1)>\dots>w(n).
\]
This rearrangement does not affect the monogressivity of $C$ (nor, for that
matter, its $ij$-ef\-fect\-ive\-ness). Indeed, the only nonobvious point here
is the relative ordering of the $a_k$ and the $b_k$: by monogressivity, we
have $w\lessdot s_{\al_k}w\lessdot s_{\al_k}s_{\al_{k'}}w$ and $w\lessdot
s_{\al_{k'}}w\lessdot s_{\al_k}s_{\al_{k'}}w$, so if $k<k'$, then, whatever
the order in which $a_k,a_{k'},b_k,b_{k'}$ appear in the original array
$[w(1)\dots w(n)]$, we must have either $a_k<b_k<a_{k'}<b_{k'}$, or
$a_{k'}<a_k<b_k<b_{k'}$. In both cases, $a_{k'}$ and $b_{k'}$ are outside of
the (numerical) interval $[a_k,b_k]$, so they may indeed appear between $a_k$
and $b_k$ in the new array without affecting monogressivity.

An orthocell thus modified will be called \emph{$ij$-normal}. The proof will
be finished if we show that there are $D_{n;i,j}$ $ij$-normal orthocells in
$S_n$. Since we have the recursion rule
\[
D_{n+1;i,j}=D_{n;i-1,j-1}+D_{n;i-1,j}+D_{n;i,j-1}+D_{n;i,j},
\]
it is enough to show that the number of $ij$-normal orthocells in $S_{n+1}$
satisfies the same recursion rule. If $C$ is such an orthocell, there are two
possibilities.
\begin{itemize}
\item Either each $s_{\al_k}$ fixes $n+1$. Removing $n+1$ from the array
$[w(1)\dots w(n+1)]$, we then obtain an orthocell in $S_n$, which is
$(i-1)(j-1)$-normal, $i(j-1)$-normal, or $ij$-normal, according to the
position of $n+1$ in the array, relative to $w(i)$ and $w(j)$.
\item Or some $s_{\al_k}$ involves $n+1$: necessarily, $k=d$ and $b_d=n+1$.
Removing again $n+1$ from the array, and discarding $\al_d$ from $C$,
we then obtain an $(i-1)j$-normal orthocell in $S_n$.
\end{itemize}
Clearly, this procedure may be reversed, starting from a normal orthocell in
$S_n$ and inserting $n+1$ at all possible places in the corresponding array.
Hence the desired recursion rule.
\end{proof}

Now let us recall a presentation for the quantized enveloping algebra
$\Uqsln$, following e.g.~\cite{Ja}: it is generated by $4(n-1)$ elements
$K_\be$, $K_\be\inv$, $X_\be$, $Y_\be$ ($\be$ a simple root), subject to the
commutation relations
\begin{align*}
K_\be K_\be\inv=1&=K_\be\inv K_\be,\qquad K_\be K_\ga=K_\ga K_\be,
\\
K_\be X_\ga K_\be\inv&=q^{(\be|\ga)} X_\ga,
\\
K_\be Y_\ga K_\be\inv&=q^{-(\be|\ga)} Y_\ga,
\\
X_\be Y_\ga-Y_\ga X_\be&=\de_{\be\ga}\,\frac{K_\be-K_\be\inv}{q-q\inv},
\end{align*}
as well as the quantized Serre relations
\begin{align*}
X_\be^2X_\ga-(q+q\inv)X_\be X_\ga X_\be+X_\ga X_\be^2&=0
\qquad\text{if $\be,\ga$ adjacent,}
\\
X_\be X_\ga-X_\ga X_\be&=0
\qquad\text{if $\be,\ga$ not adjacent,}
\\[2mm]
Y_\be^2Y_\ga-(q+q\inv)Y_\be Y_\ga Y_\be+Y_\ga Y_\be^2&=0
\qquad\text{if $\be,\ga$ adjacent,}
\\
Y_\be Y_\ga-Y_\ga Y_\be&=0
\qquad\text{if $\be,\ga$ not adjacent.}
\end{align*}
Moreover, $\Uqsln$ is a Hopf algebra whose comultiplication is given on the
generators by
\begin{equation}
\label{comult}
\begin{aligned}
\D K_\be^{\pm1}&=K_\be^{\pm1}\o K_\be^{\pm1},
\\
\D X_\be&=X_\be\o 1+K_\be\o X_\be,
\\
\D Y_\be&=Y_\be\o K_\be\inv+1\o Y_\be.
\end{aligned}
\end{equation}
We then define a $\Uqsln$-mod\-ule structure on $V^i$ as follows. For every
$w\in W$ and every simple root $\be$, we set
\[
K_\be e_w^i=q^{(w\om_i|\be)}e_w^i,
\qquad
K_\be\inv e_w^i=q^{-(w\om_i|\be)}e_w^i,
\]
and
\begin{align*}
X_\be e_w^i&=0,&Y_\be e_w^i&=e_{s_\be w}^i&&\text{if
$(w\om_i|\be)=1$;}
\\
X_\be e_w^i&=0,&Y_\be e_w^i&=0&&\text{if $(w\om_i|\be)=0$;}
\\
X_\be e_w^i&=e_{s_\be w}^i,&Y_\be e_w^i&=0&&\text{if
$(w\om_i|\be)=-1$.}
\end{align*}
(These are the only possible values for $(w\om_i|\be)$, because $\om_i$ is
minuscule.) It is straightforward to check that this module structure is well
defined, and that it is the simple $\Uqsln$-mod\-ule of highest weight
$\om_i$.
\begin{lem}
\label{ijinv}
The subspace $V^{ij}$ is a $\Uqsln$-sub\-mod\-ule of $V^i\o V^j$.
\end{lem}
\noindent By Lemma~\ref{spanned}, the statement means that the action of a
generator of $\Uqsln$ on a vector $e_C^{ij}$, $C$ monogressive and
$ij$-ef\-fect\-ive, must again be a linear combination of such vectors. We
postpone these rather tedious computations to Appendix~\ref{ijinvproof}.
\begin{cor}
\label{hwij}
The subspace $V^{ij}$ is equal to the (unique) $\Uqsln$-sub\-mod\-ule of
$V^i\o V^j$ of highest weight $\om_i+\om_j$, and the $e_C^{ij}$ (for $C$
monogressive and $ij$-ef\-fect\-ive) are linearly independent.
\end{cor}
\begin{proof}
The vector $e_{C(\rien;1)}^{ij}=e_1^i\o e_1^j$ is a highest weight
vector, of weight $\om_i+\om_j$. Now apply Lemmas \ref{spanned}~and
\ref{rightnumber}, noting that the dimension of the simple module of highest
weight $\om_i+\om_j$ is precisely $D_{n;i,j}$.
\end{proof}
\begin{lem}
\label{flipmodiso}
The linear map $R^{ji}:V^{ji}\to V^{ij}$ defined by
\begin{equation}
\label{Rijdef}
R^{ji}(e_C^{ji})=e_C^{ij}\qquad\text{for all monogressive $ij$-ef\-fect\-ive
$C$}
\end{equation}
is an isomorphism of $\Uqsln$-mod\-ules.
\end{lem}
\begin{proof}
This is immediate from the action of the generators of $\Uqsln$ on the basis
elements of $V^{ij}$ and $V^{ji}$, as described in Appendix~\ref{ijinvproof}:
the formulas obtained there are symmetric in $i$ and $j$.
\end{proof}
Extend $R^{ji}$ to an isomorphism $V^j\o V^i\simto V^i\o V^j$ of
$\Uqsln$-mod\-ules (in an arbitrary way).
\begin{lem}
\label{flipmatch}
The maps $R^{ji}$ induce isomorphisms
$R_{ij}:\L_i\o\L_j^{\sig_i}\simto\L_j\o\L_i^{\sig_j}$ of line bundles over
$E^\DrJ$, and the latter satisfy \eqref{Lijkbraid} for all $i,j,k$.
\end{lem}
\begin{proof}
The first statement amounts to the commutativity of the
diagram~\eqref{crossdiag}, which immediately follows from \eqref{imagep}~and
\eqref{Rijdef}.

For the second statement, consider the composite map
\begin{multline}
E^\DrJ\xto{\diag}E^\DrJ\times E^\DrJ\times E^\DrJ
\xto{\id\times\sig_i\times\sig_i\sig_j}
E^\DrJ\times E^\DrJ\times E^\DrJ
\\
\xto{\Pl^i\times\Pl^j\times\Pl^k}
\P(V^i)\times\P(V^j)\times\P(V^k)
\xto{\Segre}\P(V^i\o V^j\o V^k)
\end{multline}
corresponding to the line bundle $\L_i\o\L_j^{\sig_i}\o\L_k^{\sig_i\sig_j}$,
and denote by $V^{ijk}\subset V^i\o V^j\o V^k$ the linear span of the image
of this map.

{\bf Claim~A.} \emph{The subspace $V^{ijk}$ is contained in the unique simple
$\Uqsln$-sub\-module $W^{ijk}$ of $V^i\o V^j\o V^k$ of highest weight
$\om_i+\om_j+\om_k$.}

Indeed, let $K_{ij}$ be the kernel of a projection $V_i\o V_j\to
V_{\om_i+\om_j}$, and define similarly $K_{jk}$, $K_{ijk}$. Since $M^\DrJ$ is
quadratic (cf.~\cite{TT}), we have $K_{ijk}=K_{ij}\o V_k+V_i\o K_{jk}$, so
dually, $W^{ijk}=V^{ij}\o V^k\cap V^i\o V^{jk}$, and $V^{ijk}$ is clearly
contained in the right hand side. This shows Claim~A.

The proof will be finished if we show the following Claim (from which
\eqref{Lijkbraid} follows):

{\bf Claim~B.} \emph{Consider the maps $(R^{ji}\o\id)(\id\o
R^{ki})(R^{kj}\o\id)$ and $(\id\o R^{kj})(R^{ki}\o\id)(\id\o R^{ji})$ from
$V^k\o V^j\o V^i$ to $V^i\o V^j\o V^k$. Their restrictions to $V^{kji}$
agree.}

By Claim~A, it will be enough to show that the restrictions to $W^{kji}$
agree. Since both maps are morphisms between the simple $\Uqsln$-mod\-ules
$W^{kji}$ and $W^{ijk}$, they must be equal up to a constant. But they both
send the (highest weight) vector $e_1^k\o e_1^j\o e_1^i$ to $e_1^i\o e_1^j\o
e_1^k$, so this constant is equal to~$1$. This shows Claim~B.
\end{proof}
It now follows from Lemma~\ref{flipmatch} that the tuple
$T^\DrJ=(E^\DrJ,\sig_1,\dots,\sig_\ell,\L_1,\dots,\L_\ell)$ is braided.
It also follows from Corollary~\ref{hwij}~and Lemma~\ref{flipmodiso} that the
quadratic algebras $M(T^\DrJ)$ and $M^\DrJ$ agree (as quotients of
$T(V_1\oplus\dots\oplus V_{n-1})$). Conjecture~\ref{mainconj} is thus proved
for $G=\SL(n)$.

\appendix

\section{Proof of Proposition~\ref{assocbraid}}
\label{assocbraidproof}

Assume that $M_A$ is associative. By Definition~\ref{Gqdef}(c), there exists
an $A$-iso\-morph\-ism $R_{ij}:V_i\o V_j\to V_j\o V_i$. Rescaling $R_{ij}$ if
necessary, we may assume that the diagram
\begin{equation}
\label{scaleRij}
\xymatrix@C=5mm{
V_i\o V_j\ar[dr]_m\ar[rr]^{R_{ij}}&&V_j\o V_i\ar[dl]^m
\\
&V_{\om_i+\om_j}
}
\end{equation}
commutes. Now consider the following diagram:
\[
\xymatrix@C=-4mm{
&&&
V_i V_j V_k\ar[drrr]\ar[dr]\ar[dl]\ar[dlll]
\\
V_j V_i V_k\ar[rr]\ar[dr]\ar[dd]
&&
V_{i+j} V_k\ar[dr]
&&
V_i V_{j+k}\ar[dl]
&&
V_i V_k V_j\ar[ll]\ar[dl]\ar[dd]
\\
&
V_j V_{i+k}\ar[rr]
&&
V_{i+j+k}
&&
V_{i+k} V_j\ar[ll]
\\
V_j V_k V_i\ar[ur]\ar[rr]\ar[drrr]
&&
V_{j+k} V_i\ar[ur]
&&
V_k V_{i+j}\ar[ul]
&&
V_k V_i V_j\ar[ul]\ar[ll]\ar[dlll]
\\
&&&
V_k V_j V_i,\ar[ul]\ar[ur]
}
\]
where we have omitted all tensor product symbols and written $V_{i+j}$
instead of $V_{\om_i+\om_j}$, etc. (The arrows are the obvious ones, coming
either from the multiplication $m$ or from the $R_{ij}$.) All diamonds
commute by associativity, and all triangles commute, being instances
of \eqref{scaleRij}. Moreover, each object in the diagram contains a unique
copy of $V_{\om_i+\om_j+\om_k}$, and when all arrows are restricted to these
subcomodules, they become isomorphisms. Therefore, the outer rim commutes,
i.e.\ \eqref{Vijkbraid} holds.

Conversely, assume that \eqref{Vijkbraid} holds for all $i>j>k$. We first
extend the definition of the $R_{ij}$ by setting $R_{ii}:=\id$ for all $i$
and $R_{ji}:=R_{ij}\inv$ for all $i>j$. Relation~\eqref{Vijkbraid}
then holds for all $i,j,k$.

We will realize $M_A$ as a quotient of the tensor
algebra $\MT$. Let $\Gamma$ be the free monoid on
$\{1,\dots,\ell\}$. For every $I=i_1\dots i_r\in\Gamma$, let
\[
\om_I:=\om_{i_1}+\dots+\om_{i_r},\qquad
V^{\o I}:=V_{i_1}\o\dots\o V_{i_r},
\]
so $\MT=\bigoplus_{I\in\Gamma} V_I$. By Definition~\ref{Gqdef}(c), $V^{\o I}$
contains a unique copy of $V_{\om_I}$; let $K_I\subset V^{\o I}$ be its
unique invariant supplement. Then the direct sum $K=\bigoplus_{I\in\Gamma}K_I$
is a (two-sided) ideal in $\MT$ (again by Definition~\ref{Gqdef}(c)), so we
get a quotient algebra
\[
\MT/K=:\bigoplus_{I\in\Gamma}V_I.
\]
We still need to identify $V_I$ with $V_J$ whenever $\om_I=\om_J$. The
argument will be fairly standard: use the $R_{ij}$ to exchange generators from
different $V_i$'s, and check that this is consistent,
using~\eqref{Vijkbraid}. More explicitly, denote by $S_r$ the symmetric group.
If $\pi\in S_r$ and $I=i_1\dots i_r\in\Gamma$, define $\pi I:=i_{\pi(1)}\dots
i_{\pi(r)}$. Denote the usual generators of $S_r$ by $s_j:=(j,j+1)$, $1\le
j\le r-1$, and define an $A$-iso\-morph\-ism $R_{I,j}:V^{\o I}\to V^{\o
s_jI}$ by $R_{i_ji_{j+1}}$ on $V_{i_j}\o V_{i_{j+1}}$ and by $\id$ on all
other $V_{i_k}$. If $\pi\in S_r$ decomposes as $\pi=s_{j_1}\dots s_{j_t}$
(\emph{not} necessarily in a reduced way), define $R_{I,\pi}:V^{\o I}\to
V^{\o\pi I}$ by $R_{I,\pi}:=R_{I,j_1}\dots R_{I,j_t}$. Since this is an
$A$-iso\-morph\-ism, it restricts to $R_{I,\pi}:V_I\to V_{\pi I}$. Note that
this restriction does not depend on the chosen decomposition of $\pi$, thanks
to $R_{ii}=\id$, to $R_{ij}R_{ji}=\id$, and to \eqref{Vijkbraid}.

Now consider $\MT/K$ as a $\Pp$-graded algebra, the term of degree
$\la\in\Pp$ being $U_\la:=\bigoplus_{I:\om_I=\la}V_I$. Let $S_\la\subset
U_\la$ be the span of all elements $x-R_{I,\pi}(x)$, $x\in V_I$, where $I$
runs over all elements of $\Gamma$ such that $\om_I=\la$. Then $U_\la/S_\la$
consists of just one copy of $V_\la$: indeed, on one hand, $\om_I=\om_J$ if
and only if $J=\pi I$ for some $\pi\in S_r$, and on the other hand, the
construction of the $R_{I,\pi}$ implies that
\[
R_{I,\pi\pi'}=R_{\pi'I,\pi}R_{I,\pi'}.
\]
Moreover, the direct sum $S:=\bigoplus_{\la\in\Pp}S_\la$ is an ideal in
$\MT/K$, so the corresponding quotient yields the desired associative
realization of the shape algebra $M_A$.

\section{Proof of Proposition~\ref{glue}}
\label{glueproof}

We need to show that for any monogressive orthocells $C_1,C_2$, the
automorphisms $\sig_{i,C_1}$ and $\sig_{i,C_2}$ agree on $E(C_1)\cap
E(C_2)$. Since this intersection is clearly $T$-stable and closed, it is a
union of $T$-orbit closures in each of $E(C_1)$ and $E(C_2)$,
hence~\cite[Corollary~6.2]{PPinGB} a union of $E(C')$ with $C'$ a common
subcell of $C_1$ and $C_2$.

It is therefore enough to show that the $\sig_{i,C}$ are compatible with
restriction to subcells. Consider a monogressive orthocell
$C=C(w;\al_1,\dots,\al_d)$ and a subcell, say,
$C'=C(ws_L;\al'_1,\dots,\al'_e)$, with $L\subset\{1,\dots,d\}$,
$\{\al'_1,\dots,\al'_e\}\subset\{\al_1,\dots,\al_d\}$, and
$\al_k\not\in\{\al'_1,\dots,\al'_e\}$
for all $k\in L$ (so that $ws_L$ is again of minimal length in $C'$). For
each $\al\in\Phi$, fix an isomorphism $u_\al:(\C,{+})\to U_\al$ such that
$tu_\al(z)t\inv=u_\al(\al(t)z)$ for all $t\in T$ and all
$z\in\C$~\cite[Proposition~8.1.1(i)]{Sp}. Then the
set of all $\dot w\dot s_L u_{-\al'_1}(z_1)\dots u_{-\al'_e}(z_e)B$,
$(z_1,\dots,z_e)\in\C^e$, is an open dense subset of $E(C')$
(cf.~\cite{PPinGB}, proof of Theorem~4.1), and the action of $\sig_{i,C'}$
on such an element is given by
\begin{align*}
&ws_Lt_i s_L\inv w\inv\dot w\dot s_L u_{-\al_1}(z_1)\dots u_{-\al_d}(z_d)B
\\
&\qquad=\dot w\dot s_L t_i u_{-\al_1}(z_1)\dots u_{-\al_d}(z_d)B
\\
&\qquad=\dot w\dot s_L
u_{-\al_1}(\al_1(t_i)\inv z_1)\dots u_{-\al_d}(\al_d(t_i)\inv z_d)B,
\end{align*}
whereas the action of $\sig_{i,C}$ (i.e.\ multiplication by $wt_i w\inv$
instead of $ws_Lt_i s_L\inv w\inv$) is given by the same expression, with
$t_i$ replaced by $s_L\inv t_i s_L$. But since $s_L$ is a product of
reflections w.r.t.\ roots orthogonal to each $\al'_k$, we have
$\al'_k(s_L\inv ts_L)=(s_L\al'_k)(t)=\al'_k(t)$. Thus, the restriction of
$\sig_{i,C}$ to $C'$ coincides with $\sig_{i,C'}$, and the result follows.

\section{Proof of Lemma~\ref{ijinv}}
\label{ijinvproof}

We begin by collecting some more explicit information on the root system of
$\SL(n)$. First, $(\al|\al)=2$ for every root $\al$, so in particular,
\[
s_\al(\la)=\la-(\la|\al)\,\al
\]
for any weight $\la$. Recall also that all fundamental weights
$\om_1,\dots,\om_{n-1}$ are minuscule, so for any $w\in S_n$ and any root
$\al$, $(w\om_i|\al)=0$, $1$, or $-1$. Moreover, if $\al>0$, then
\begin{align*}
w<s_\al w&\implies\text{$(w\om_i|\al)=0$ or $1$,}
\\
w>s_\al w&\implies\text{$(w\om_i|\al)=0$ or $-1$.}
\end{align*}
Now let $\al\ne\al'$ be two positive roots and $s_\al=(a\,b)$,
$s_{\al'}=(a'\,b')$, with $a<b$ and $a'<b'$. Then
\[
(\al|\al')=
\begin{cases}
1&\text{if $a=a'$ or $b=b'$ (but not both),}
\\
0&\text{if $\{a,b\}\cap\{a',b'\}=\emptyset$,}
\\
-1&\text{if $a=b'$ or $b=a'$.}
\end{cases}
\]
The preceding information will be used freely in the sequel, without explicit
reference.

We fix a simple root $\be$. Consider first the action of the generator
$K_\be$ on a vector $e_C^{ij}$, where $C=C(\al_1,\dots,\al_d;w)$ is
monogressive and $ij$-ef\-fect\-ive. Recalling the expression \eqref{comult}
for $\D K_\be$, we get
\[
K_\be e_C^{ij}=\sum_{L\subset\{1,\dots,d\}}
q^{(s_{\bar L}w\om_i|\be)+(s_Lw\om_j|\be)} q^{|L|}\,e_{s_{\bar L}w}^i\o
e_{s_Lw}^j.
\]
For each $1\le k\le d$, we have $s_{\al_k}w\om_j=w\om_j-\al_k$. More
generally, $s_Lw\om_j=w\om_j-\sum_{k\in L}\al_k$, and similarly for $s_{\bar
L}w\om_i$; therefore,
\[
K_\be e_C^{ij}=q^{(w(\om_i+\om_j)|\be)-\sum_{k=1}^d(\al_k|\be)}\,e_C^{ij}.
\]
A similar formula holds for $K_\be\inv e_C^{ij}$.

Now we study the action of $X_\be$ and of $Y_\be$ on a vector $e_C^{ij}$.
Note that the root $\be$ will be orthogonal to all defining roots of the
orthocell $C$, except at most two. Thus, there are four cases to consider:
\begin{description}
\item[Case I] $C=C(\al,\al',\al_1,\dots,\al_d;w)$,
\item[Case II] $C=C(\al,\al_1,\dots,\al_d;w)$,
\item[Case III] $C=C(\al_1,\dots,\al_d;w)$,
\item[Case IV] $C=C(\be,\al_1,\dots,\al_d;w)$,
\end{description}
where, in all cases, $\al,\al',\al_1,\dots,\al_d$ are pairwise orthogonal,
$(\be|\al)=\pm1$, $(\be|\al')=\pm1$, and $(\be|\al_k)=0$ for all $k$.

We will first treat these four cases when $d=0$, and then describe how to
deduce results for arbitrary $d$ from this particular case.

Let us use the notation $c\lessdot c'$ even when $c,c'$ are integers, meaning
that $c'=c+1$. We will also use the following notation throughout:
\[
s:=s_\al=:(a\,b),
\qquad
s':=s_{\al'}=:(a'\,b'),
\qquad
t:=s_\be,
\]
with $a<b$ and $a'<b'$. Note that monogressivity and $ij$-ef\-fect\-ive\-ness
\emph{exclude} the orderings $a<a'<b<b'$ and $a'<a<b'<b$.
\vspace{5mm}

\noindent{\bf Case I: $C=C(\al,\al';w)$.}
\nopagebreak
\vspace*{5mm}

\noindent\emph{Subcase I.1: $(\be|\al)=(\be|\al')=-1$:}
Exchanging $\al,\al'$ if necessary, we may assume that $a<a'$. Then we must
have $a<b\lessdot a'<b'$ and $t=(b\,a')$. Furthermore,
$(sw\om_i|\be)=(s'w\om_i|\be)=(w\om_i|\be)-1$ and
$(ss'w\om_i|\be)=(w\om_i|\be)-2$. Since all these inner products
must be equal to $0$, $1$, or $-1$, we get
\[
(w\om_i|\be)=-1,
\qquad(sw\om_i|\be)=(s'w\om_i|\be)=0,
\qquad(ss'w\om_i|\be)=1,
\]
and similarly for $\om_i$ replaced by $\om_j$. Recalling the expression
\eqref{comult} for $\D X_\be$ and $\D Y_\be$, we obtain
\begin{align*}
X_\be e_C^{ij}
&=X_\be\left(q^2e_w^i\o e_{ss'w}^j
+q\,e_{sw}^i\o e_{s'w}^j
+q\,e_{s'w}^i\o e_{sw}^j
+e_{ss'w}^i\o e_w^j\right)
\\
&=q^2X_\be e_w^i\o e_{ss'w}^j
+K_\be e_{ss'w}^i\o X_\be e_w^j
\\
&=q^2e_{tw}^i\o e_{ss'w}^j
+q\,e_{ss'w}^i\o e_{tw}^j,
\end{align*}
and, by a similar computation,
\[
Y_\be e_C^{ij}
=q^2e_w^i\o e_{tss'w}^j
+q\,e_{tss'w}^i\o e_w^j.
\]
The vanishing inner products obtained above imply that $tsw\om_i=sw\om_i$ and
$ts'w\om_i=s'w\om_i$. Using the Coxeter relations $(ts)^3=(ts')^3=1$,
we then also have the equalities
\begin{align*}
tw\om_i&=stw\om_i=s'tw\om_i=ss'tw\om_i,
\\
ss'w\om_i&=sts'w\om_i=s'tsw\om_i=ss'tss'tw\om_i,
\\
w\om_i&=tstw\om_i=ts'tw\om_i=tss'tw\om_i,
\\
tss'w\om_i&=stss'w\om_i=s'tss'w\om_i=ss'tss'w\om_i
\end{align*}
(and similarly for $\om_i$ replaced by $\om_j$), which may be used in the
above expressions for $X_\be e_C^{ij}$ and $Y_\be e_C^{ij}$; cf.\
Remark~\ref{equiveff}. Since
$a<b<a'<b'$, monogressivity implies that there are four possible relative
positions of $a,b,a',b'$ inside the array $w$, yielding the following
expressions for $X_\be e_C^{ij}$ and $Y_\be e_C^{ij}$:
\begin{center}
\begin{tabular}{c|c|c}
$w$&$X_\be e_C^{ij}$&$Y_\be e_C^{ij}$:
\\
\hline
$[\cdot a\cdot a'\cdot b\cdot b'\cdot]$
&
$-q\,e_{C(ss'(\be);ss'tw)}^{ij}$
&
$-q\,e_{C(ss'(\be);tss'tw)}^{ij}$,
\\
$[\cdot a'\cdot a\cdot b'\cdot b\cdot]$
&
$-q\,e_{C(ss'(\be);tw)}^{ij}$
&
$-q\,e_{C(ss'(\be);w)}^{ij}$,
\\
$[\cdot a'\cdot a\cdot b\cdot b'\cdot]$
&
$-q\,e_{C(ss'(\be);s'tw)}^{ij}$
&
$-q\,e_{C(ss'(\be);ts'tw)}^{ij}$,
\\
$[\cdot a\cdot a'\cdot b'\cdot b\cdot]$
&
$-q\,e_{C(ss'(\be);stw)}^{ij}$
&
$-q\,e_{C(ss'(\be);tstw)}^{ij}$.
\end{tabular}
\end{center}
These results are valid provided all orthocells involved are monogressive
(their $ij$-ef\-fect\-ive\-ness being clear). In each case, this may easily
be checked using Criterion~\ref{monocrit}. For example, in the first line, we
have $w=[\cdot a\cdot a'\cdot b\cdot b'\cdot]$, $ss'tw=[\cdot b\cdot a\cdot
b'\cdot a'\cdot]$, and $s_{ss'(\be)}=(a\,b')$: since $C=C(\al,\al';w)$ is
monogressive by assumption, the subarray $[a'\dots b]$ of $w$ contains no
numbers in the (numerical) interval $[a,b]$, nor in $[a',b']$, and therefore
not in $[a,b']$ (because $a<b\lessdot a'<b'$); thus, the orthocell
$C(ss'(\be);ss'tw)$ is again monogressive.
\vspace{5mm}

\noindent\emph{Subcase I.2: $(\be|\al)=1$ and $(\be|\al')=-1$:}
Here we must have either $a\lessdot a'<b'<b$ and $t=(a\,a')$, or
$a<a'<b'\lessdot b$ and $t=(b'\,b)$. Furthermore,
\[
(w\om_i|\be)=(ss'w\om_i|\be)=0,
\qquad(sw\om_i|\be)=-1,
\qquad(s'w\om_i|\be)=1,
\]
and similarly for $\om_i$ replaced by $\om_j$. It follows that
\begin{align*}
X_\be e_C^{ij}&=q\,e_{tsw}^i\o e_{s'w}^j
+q^2e_{s'w}^i\o e_{tsw}^j,
\\
Y_\be e_C^{ij}&=q\,e_{sw}^i\o e_{ts'w}^j
+q^2e_{ts'w}^i\o e_{sw}^j.
\end{align*}
Arguments similar to those of Subcase~I.1 then show (whether $t=(a\,a')$ or
$t=(b'\,b)$) that the orthocells $C(ss'(\be);s'tw)$ and $C(ss'(\be);ts'tw)$
are monogressive and that
\[
X_\be e_C^{ij}=-q\,e_{C(ss'(\be);s'tw)}^{ij},
\qquad
Y_\be e_C^{ij}=-q\,e_{C(ss'(\be);ts'tw)}^{ij}.
\]
(We omit the details.)
\vspace{5mm}

\noindent\emph{Subcase I.3: $(\be|\al)=(\be|\al')=1$:}
These inner products force $a<a'<b<b'$ or $a'<a<b'<b$, contradicting
the fact that $C(\al,\al';w)$ is monogressive and $ij$-ef\-fect\-ive.
Therefore, this Subcase is impossible.
\vspace{5mm}

\noindent{\bf Case II: $C=C(\al;w)$.}
\nopagebreak
\vspace*{5mm}

\noindent\emph{Subcase II.1: $(\be|\al)=-1$:}
We must have either $c\lessdot a<b$ and $t=(c\,a)$, or  $a<b\lessdot c$ and
$t=(b\,c)$. Moreover, since $(sw\om_i|\be)=(w\om_i|\be)+1$, and similarly for
$\om_j$, we obtain the following cases.
\begin{itemize}
\item If $(w\om_i|\be)=(w\om_j|\be)=-1$, then arguments similar to those of
Case~I show that
\[
Y_\be e_C^{ij}=0
\]
and that $X_\be e_C^{ij}=q\,e_{tw}^i\o e_{sw}^j+e_{sw}^i\o e_{tw}^j$ is given
by the following table:
\begin{center}
\begin{tabular}{c|c|c}
$t$&$w$&$X_\be e_C^{ij}$:
\\
\hline
$(c\,a)$&$[\cdot a\cdot b\cdot c\cdot]$&$e_{C(s(\be);tw)}^{ij}$,
\\
$(c\,a)$&$[\cdot a\cdot c\cdot b\cdot]$&$e_{C(s(\be);stw)}^{ij}$,
\\
$(b\,c)$&$[\cdot a\cdot c\cdot b\cdot]$&$-e_{C(s(\be);stw)}^{ij}$,
\\
$(b\,c)$&$[\cdot c\cdot a\cdot b\cdot]$&$-e_{C(s(\be);tw)}^{ij}$.
\end{tabular}
\end{center}
\item If $(w\om_i|\be)=(w\om_j|\be)=0$, then
\[
X_\be e_C^{ij}=0
\]
and $Y_\be e_C^{ij}=q\,e_w^i\o e_{tsw}^j+e_{tsw}^i\o e_w^j$ is given by
the following table:
\begin{center}
\begin{tabular}{c|c|c}
$t$&$w$&$Y_\be e_C^{ij}$:
\\
\hline
$(c\,a)$&$[\cdot c\cdot a\cdot b\cdot]$&$-e_{C(s(\be);tw)}^{ij}$,
\\
$(c\,a)$&$[\cdot a\cdot c\cdot b\cdot]$&$-e_{C(s(\be);w)}^{ij}$,
\\
$(b\,c)$&$[\cdot a\cdot c\cdot b\cdot]$&$e_{C(s(\be);w)}^{ij}$,
\\
$(b\,c)$&$[\cdot a\cdot b\cdot c\cdot]$&$e_{C(s(\be);tw)}^{ij}$.
\end{tabular}
\end{center}
\item If $(w\om_i|\be)=-1$ and $(w\om_j|\be)=0$, or vice versa, then
\[
X_\be e_C^{ij}=-q\,e_{C(\rien;stw)}^{ij},
\qquad
Y_\be e_C^{ij}=-q\,e_{C(\rien;stsw)}^{ij}.
\]
\end{itemize}
\vspace{5mm}

\noindent\emph{Subcase II.2: $(\be|\al)=1$:}
Here we have either $a\lessdot c<b$ and $t=(a\,c)$, or $a<c\lessdot b$ and
$t=(c\,b)$. Moreover, since $(sw\om_i|\be)=(w\om_i|\be)-1$, and similarly for
$\om_j$, we obtain the following cases.
\begin{itemize}
\item If $(w\om_i|\be)=(w\om_j|\be)=1$, then 
\[
X_\be e_C^{ij}=0,
\qquad
Y_\be e_C^{ij}=\pm e_{C(s(\be);tw)}^{ij}.
\]
\item If $(w\om_i|\be)=(w\om_j|\be)=0$, then
\[
X_\be e_C^{ij}=\pm e_{C(s(\be);tw)}^{ij},
\qquad
Y_\be e_C^{ij}=0.
\]
\item The case $(w\om_i|\be)=1$ and $(w\om_j|\be)=0$, or vice versa,
contradicts the $ij$-ef\-fect\-ive\-ness of $C$ and is therefore impossible.
\end{itemize}
\vspace{5mm}

\noindent{\bf Case III: $C=C(\rien;w)$.}
We obtain the following table:
\begin{center}
\begin{tabular}{c|c|c|c}
$(w\om_i|\be)$&$(w\om_j|\be)$&$X_\be e_C^{ij}$&$Y_\be e_C^{ij}$:
\\
\hline
$1$ & $1$ & $0$ & $q\inv\,e_{C(\be;w)}^{ij}$,
\\
$1$ & $0$ & $0$ & $e_{C(\rien;tw)}^{ij}$,
\\
$0$ & $1$ & $0$ & $e_{C(\rien;tw)}^{ij}$,
\\
$0$ & $0$ & $0$ & $0$,
\\
$0$ & $-1$ & $e_{C(\rien;tw)}^{ij}$ & $0$,
\\
$-1$ & $0$ & $e_{C(\rien;tw)}^{ij}$ & $0$,
\\
$-1$ & $-1$ & $q\inv\,e_{C(\be;tw)}^{ij}$ & $0$.
\end{tabular}
\end{center}
\vspace{5mm}

\noindent{\bf Case IV: $C=C(\be;w)$.}
Since $C$ is monogressive and $ij$-ef\-fect\-ive, we must have
$(w\om_i|\be)=(w\om_j|\be)=1$, hence
\[
X_\be e_C^{ij}=(q^2+1)e_{C(\rien;w)}^{ij},
\qquad
Y_\be e_C^{ij}=(q^2+1)e_{C(\rien;tw)}^{ij}.
\]

Finally, we show how, in the preceding four cases, one can deduce the action
of $X_\be$ and $Y_\be$ for arbitrary $d$ from that for $d=0$. The idea is
that $\al_1,\dots,\al_d$, being orthogonal to $\be,\al,\al'$, do not
``interfere'' with the computations done above. To make this idea precise, we
will restrict ourselves to the very first case treated above (all other cases
being similar), namely, the action of $X_\be$ on $e_C^{ij}$ when
$C=C(\al,\al',\al_1,\dots,\al_d;w)$, $(\be|\al)=(\be|\al')=-1$,
$(\be|\al_k)=0$ for all $k$, $s_\al=(a\,b)$ and $s_{\al'}=(a'\,b')$ with
$a<b\lessdot a'<b'$ (so $t:=s_\be=(b\,a')$), and $w=[\cdot a\cdot a'\cdot
b\cdot b'\cdot]$.

Since $\be$ is orthogonal to each $\al_k$, we have $(s_L\la|\be)=(\la|\be)$
for any weight $\la$ and any $L\subset\{1,\dots,d\}$, so we still get
\[
(s_Lw\om_i|\be)=-1,
\qquad(ss_Lw\om_i|\be)=(s's_Lw\om_i|\be)=0,
\qquad(ss's_Lw\om_i|\be)=1.
\]
Therefore, the action of $X_\be$ on each term of
\begin{align*}
e_{C(\al,\al',\al_1,\dots,\al_d;w)}^{ij}
=\sum_{L\subset\{1,\dots,d\}}q^{|L|}
&\Bigl(
q^2e_{s_{\bar L}w}^i\o e_{ss's_Lw}^j
+q\,e_{ss_{\bar L}w}^i\o e_{s's_Lw}^j
\\
&\qquad
+q\,e_{s's_{\bar L}w}^i\o e_{ss_Lw}^j
+e_{ss's_{\bar L}w}^i\o e_{s_Lw}^j
\Bigr)
\end{align*}
is still computed in a similar way to that on $e_{C(\al,\al';w)}^{ij}$, viz.
\[
X_\be e_{C(\al,\al',\al_1,\dots,\al_d;w)}^{ij}
=\sum_{L\subset\{1,\dots,d\}}q^{|L|}
\left(
q^2e_{ts_{\bar L}w}^i\o e_{ss's_Lw}^j
+q\,e_{ss's_{\bar L}w}^i\o e_{ts_Lw}^j
\right).
\]
It follows that
\[
X_\be e_{C(\al,\al',\al_1,\dots,\al_d;w)}^{ij}
=-q\,e_{C(ss'(\be),\al_1,\dots,\al_d;ss'tw)}^{ij},
\]
provided the orthocell $C(ss'(\be),\al_1,\dots,\al_d;ss'tw)$ is monogressive.
But it is easy to see that the analysis of the monogressivity of
$C(ss'(\be);ss'tw)$ made earlier, using Criterion~\ref{monocrit}, remains
valid if $\al_1,\dots,\al_d$ are added to the orthocells $C(\al,\al';w)$ and
$C(ss'(\be);ss'tw)$.

\end{document}